 \documentclass{article}
\usepackage{graphicx}
\usepackage{caption}
\usepackage{amssymb,amsmath,amsthm,amscd}
\usepackage{verbatim}

\usepackage{enumitem}

\usepackage{float}

\usepackage{array}

\newcommand{\mycomment}[1]{}

\newcommand\mycup[2]{\overset{#2}{\underset{{#1}}{\cup}}}
\newcommand\mysum[2]{\overset{#2}{\underset{{#1}}{\textstyle\sum}}}

\theoremstyle{plain}
\newtheorem{theorem}{Theorem}
\newtheorem{lemma}[theorem]{Lemma}
\newtheorem{corollary}[theorem]{Corollary}
\newtheorem*{corollary*}{Corollary}
\newtheorem{proposition}[theorem]{Proposition}
\newtheorem*{claim*}{Claim}

\theoremstyle{definition}
\newtheorem{remark}[theorem]{Remark} 

\theoremstyle{plain}
\makeatletter
\newtheorem*{rep@theorem}{\rep@title}
\newcommand{\newreptheorem}[2]{%
\newenvironment{rep#1}[1]{%
 \def\rep@title{#2 \ref{##1}}%
 \begin{rep@theorem}}%
 {\end{rep@theorem}}}
\makeatother

\newreptheorem{theorem}{Theorem}
\newreptheorem{lemma}{Lemma}
\newreptheorem{proposition}{Proposition}
\newreptheorem{corollary}{Corollary}

\def\R{{\mathbb R}}
\def\Z{{\mathbb Z}}

\def\cL{{\mathcal L}}
\def\U{{\mathcal{U}}}

\newcommand*{\defeq}{\mathrel{\vcenter{\baselineskip0.5ex \lineskiplimit0pt
                     \hbox{\scriptsize.}\hbox{\scriptsize.}}}%
                     =}

\usepackage{mathtools}

\makeatletter
\DeclareRobustCommand\widecheck[1]{{\mathpalette\@widecheck{#1}}}
\def\@widecheck#1#2{%
    \setbox\z@\hbox{\m@th$#1#2$}%
    \setbox\tw@\hbox{\m@th$#1%
       \widehat{%
          \vrule\@width\z@\@height\ht\z@
          \vrule\@height\z@\@width\wd\z@}$}%
    \dp\tw@-\ht\z@
    \@tempdima\ht\z@ \advance\@tempdima2\ht\tw@ \divide\@tempdima\thr@@
    \setbox\tw@\hbox{%
       \raise\@tempdima\hbox{\scalebox{1}[-1]{\lower\@tempdima\box
\tw@}}}%
    {\ooalign{\box\tw@ \cr \box\z@}}}
\makeatother

\renewcommand\sup[1]{^{({#1})}}

\newcommand{\eps}{\varepsilon}

\newcommand{\sign}{\operatorname{sign}}
\newcommand{\self}{\operatorname{self}}

\renewcommand{\int}{\operatorname{int}}
\newcommand{\lk}{\operatorname{lk}}
\renewcommand{\d}{\partial}

\newcommand{\id}{\operatorname{id}}

\renewcommand\setminus{\mathbin{\mathpalette\setminusaux\relax}}
\newcommand\setminusaux[2]{\mspace{-4mu}
  \raisebox{\rsmraise{#1}\depth}{\rotatebox[origin=c]{-20}{$#1\smallsetminus$}}
 \mspace{-4mu}
}
\newcommand\rsmraise[1]{%
  \ifx#1\displaystyle .8\else
    \ifx#1\textstyle .8\else
      \ifx#1\scriptstyle .6\else
        .45%
      \fi
    \fi
  \fi}

\usepackage{xcolor}
\newcommand*{\missingreference}[1]{\colorbox{yellow}{?#1?}}
\newcommand*{\missingcitation}[1]{\colorbox{green}{?#1?}}
\makeatletter
\def\@setref#1#2#3{%
  \ifx#1\relax
   \protect\G@refundefinedtrue
   \nfss@text{\reset@font\missingreference{#3}}
   \@latex@warning{Reference `#3' on page \thepage \space
             undefined}%
  \else
   \expandafter#2#1\null
  \fi}
\def\@citex[#1]#2{\leavevmode
  \let\@citea\@empty
  \@cite{\@for\@citeb:=#2\do
    {\@citea\def\@citea{,\penalty\@m\ }%
     \edef\@citeb{\expandafter\@firstofone\@citeb\@empty}%
     \if@filesw\immediate\write\@auxout{\string\citation{\@citeb}}\fi
     \@ifundefined{b@\@citeb}{\hbox{\reset@font\missingcitation{#2}}
       \G@refundefinedtrue
       \@latex@warning
         {Citation `\@citeb' on page \thepage \space undefined}}%
       {\@cite@ofmt{\csname b@\@citeb\endcsname}}}}{#1}}
\makeatother

\usepackage{soul}

\usepackage[caption=false,labelformat=simple]{subfig}

\mycomment{
    \usepackage{titlesec}

    \mycomment{
    \titleformat{\section}
      {\normalfont\scshape}{\thesection}{1em}{}

    \titleformat{\subsection}
      {\normalfont\scshape}{\thesubsection}{1em}{}

    \titleformat{\subsubsection}
      {\normalfont\scshape}{\thesubsubsection}{1em}{}
    }
}

\SetLabelAlign{CenterWithParen}{\hfil\makebox[1.0em]{#1}\hfil} 

\usepackage{xargs}
\usepackage[colorinlistoftodos,prependcaption,textsize=tiny]{todonotes}
\newcommandx{\unsure}[2][1=]{\todo[linecolor=red,backgroundcolor=red!25,bordercolor=red,#1]{#2}}
\newcommandx{\change}[2][1=]{\todo[linecolor=blue,backgroundcolor=blue!25,bordercolor=blue,#1]{#2}}
\newcommandx{\info}[2][1=]{\todo[linecolor=OliveGreen,backgroundcolor=OliveGreen!25,bordercolor=OliveGreen,#1]{#2}}
\newcommandx{\improvement}[2][1=]{\todo[linecolor=Plum,backgroundcolor=Plum!25,bordercolor=Plum,#1]{#2}}
\newcommandx{\thiswillnotshow}[2][1=]{\todo[disable,#1]{#2}}

\numberwithin{equation}{section} 
\numberwithin{theorem}{section}

\newcommand{\myp}{\hat} %
\newcommand{\mym}{\acute} %

\newcommand{\tightoverset}[2]{%
  \mathop{#2}\limits^{\vbox to -.5ex{\kern-0.75ex\hbox{$#1$}\vss}}}

\mycomment{
    \usepackage[parfill]

    \makeatletter
    \def\subsection{\@startsection{subsection}{3}%
      \z@{.5\linespacing\@plus.7\linespacing}{.1\linespacing}%
      {\bfseries}}
    \makeatother

    \makeatletter
    \def\subsubsection{\@startsection{subsubsection}{3}%
      \z@{.5\linespacing\@plus.7\linespacing}{.1\linespacing}%
      {\normalfont\itshape}}
    \makeatother
}

\begin{document}
\sloppy 

\title{The triviality of a certain invariant of link maps in the four-sphere}

\author{Ash Lightfoot}
\date{}  
\maketitle

\begin{abstract}
It is an open problem whether  Kirk's $\sigma$ invariant is the complete obstruction to a 
link map $S^2\cup S^2\to S^4$ being link homotopically trivial.
With the objective of constructing counterexamples, Li proposed a link homotopy invariant $\omega$ that is defined on the kernel of $\sigma$ and also obstructs link nullhomotopy. We show that $\omega$ is determined by $\sigma$, and is a strictly weaker invariant.
\end{abstract}

\section{Introduction}


A link map is a (continuous) map
\[
    f:S^{p_1}\cup S^{p_2}\cup \ldots \cup S^{p_n}\to S^m
\]
from a union of spheres into another sphere such that $f(S^{p_i})\cap f(S^{p_j})=\emptyset$ for $i\neq j$. Two link maps are said to be link homotopic if they are connected by a homotopy through link maps, and  the set of link homotopy classes of link maps as above is denoted $LM_{{p_1},p_2,\ldots, p_n}^m$. It is a familiar result that $LM_{1,1}^3$ is classified by the linking number, and in his foundational work Milnor \cite{Mi}  described invariants of $LM_{1,1,\ldots, 1}^3$ which classified $LM_{1,1,1}^3$. These invariants (the $\overline{\mu}$-invariants) were  refined much later by Habegger and Lin \cite{HL}  to achieve an algorithmic classification of $LM_{1,1,\ldots, 1}^3$.


Higher dimensional link homotopy began with a study of $LM_{p,q}^m$ when $p,q\leq m-3$, first by Scott \cite{Sc} and later by Massey and Rolfsen \cite{MR}. Both papers made particular use of a generalization of the  linking number, defined as follows. 
Given a link map $f:S^p\cup S^q\to S^m$, choose a point $\infty\in S^m\setminus f(S^p\cup S^q)$ and identify $S^m\setminus \infty$ with $\R^{m}$. When $p,q<m$, the map
\[
    S^p\times S^q\to S^{m-1}, (x,y)\mapsto \frac{f(x)-f(y)}{||f(x)-f(y)||}
\]
is nullhomotopic on the subspace $S^p\vee S^q$ and so determines an element $\alpha(f)\in \pi_{p+q}(S^{m-1})$. When $m=p+q+1$, the link homotopy invariant $\alpha$ is the integer-valued linking number.  In a certain dimension range, $\alpha$ was shown in \cite{MR}  to classify \emph{embedded} link maps $S^p\cup S^q\to S^m$ up to link homotopy. Indeed, historically,  link homotopy roughly separated into settling two problems.
\begin{enumerate}
    \item Decide when an embedded link map is link nullhomotopic.
    \item Decide when a link map is link homotopic to an embedding.
\end{enumerate}
In a large metastable range this approach culminated in a long exact sequence which reduced the problem of classifying $LM_{p,q}^m$ to standard homotopy theory questions (see \cite{Ko1}).

On the other hand, four-dimensional topology presents unique difficulties, and link homotopy of 2-spheres in the 4-sphere requires different techniques.
%
%
%
%
In this setting, the first problem listed above was solved by Bartels and Teichner \cite{BT}, who showed that an embedded link $S^2\cup S^2\cup \ldots \cup S^2\to S^4$ is link nullhomotopic. In this paper we are interested in  invariants of $LM_{2,2}^4$ which have been introduced to address the second problem.

Fenn and Rolfsen \cite{FR}  showed that $\alpha$ defines a surjection $LM_{2,2}^4\to \Z_2$ and in doing so constructed the first example of a link map $S^2\cup S^2\to S^4$ which is not link nullhomotopic. Kirk \cite{Ki1} generalized this result, introducing an invariant $\sigma$ of $LM_{2,2}^4$ which further obstructs embedding and surjects onto an infinitely generated group.

To a link map $f:S^2_+\cup S^2_-\to S^4$, where we use signs to distinguish component 2-spheres,  Kirk defined a pair of integer polynomials $\sigma(f)=(\sigma_+(f), \sigma_-(f))$ such that each component is invariant under link homotopy of $f$, determines $\alpha(f)$ and vanishes if $f$ is link homotopic to a link map that embeds \emph{either} component. It is an open problem whether $\sigma$ is the complete obstruction; that is, whether $\sigma(f)=(0,0)$ implies that $f$ is link homotopic to an embedding. By \mycomment{Since singular link concordance implies link homotopy} \cite[Theorem 5]{BT},  this is equivalent to asking if $\sigma$ is injective on $LM_{2,2}^4$.
Seeking to answer in the negative, Li proposed an invariant $\omega(f)=(\omega_+(f), \omega_-(f))$ to detect link maps in the kernel of $\sigma$.

When $\sigma_\pm(f)=0$, after a link homotopy the restricted map $f|_{S^2_{\pm}}:S^2_\pm\to S^4\setminus f(S^2_\mp)$ may be equipped with a collection of Whitney disks,  and $\omega_\pm(f)\in \Z_2$ obstructs embedding by counting weighted intersections between $f(S^2_\pm)$ and the interiors of  these disks. Precise definitions of these invariants will be given in Section \ref{sec:prelims}.
 
By \cite{Li97} (and \cite{me2}),  $\omega$ is an invariant of link homotopy, but the example produced in \cite{Li97} of a link map $f$ with $\sigma(f)=(0,0)$ and $\omega(f)\neq (0,0)$ was found to be in error by  Pilz \cite{Pilz}. The purpose of this paper is to prove that  $\omega$ cannot detect such examples;  indeed, it is a weaker invariant than $\sigma$.

\begin{theorem}\label{thm:mainresult}
Let $f$ be a link map with $\sigma_-(f)=0$ and let $a_1,a_2,\ldots$ be integers so that $\sigma_+(f)=\mysum{n\geq 1}{} a_n (s^n - 1)$.  Then
\[
    \omega_-(f) = \textstyle \sum a_n \hskip -0.2cm\mod 2,
\]
where the sum is over all $n$ equal to $2$ modulo $4$.
\end{theorem}

Consequently,  there are infinitely many distinct classes $f\in LM_{2,2}^4$ with $\sigma_+(f)=0$, $\omega_+(f)=0$ but $\sigma_-(f)\neq 0$ (see Proposition \ref{eq:PK-image}). In particular, the following corollary answers Question 6.2 of \cite{Li97}.

\mycomment{
    Maybe later: use kappa, so can say kernel infinitely generated.
}

\begin{corollary}\label{coro:maincoro}
If a link map $f$ has $\sigma(f)=(0,0)$, then $\omega(f)=(0,0)$.
\end{corollary}

By \cite[Theorem 1.3]{me2} and \cite[Theorem 2]{ST},  Theorem \ref{thm:mainresult} may be interpreted geometrically as follows.

\begin{corollary}
Let $f$ be a link map such that $\sigma_+(f)=0$. Then, after a link homotopy, the self-intersections of $f(S^2_+)$ may be paired up with framed, immersed Whitney disks in $S^4\setminus f(S^2_-)$ whose interiors are disjoint from $f(S^2_+)$.
\end{corollary}

This paper is organized as follows. In Section \ref{sec:prelims} we first review Wall intersection theory  in the four-dimensional setting. The geometric principles thus established underly the link homotopy invariants $\sigma$ and $\omega$, which we subsequently define.
In Section \ref{sec:proof} we exploit that, up to link homotopy, one component of a link map is unknotted, immersed, to equip this component with a convenient collection of Whitney disks which enable us to relate the invariants $\sigma$ and $\omega$. A more detailed outline of the proof may be found at the beginning of that section.

\section{Preliminaries}\label{sec:prelims}

Let us first fix some notation. For an oriented path or loop $\alpha$, let $\overline{\alpha}$ denote its reverse path; if $\alpha$ is a based loop, let $[\alpha]$ denote its based homotopy class. Let $\ast$ denote composition of paths, and denote the interval $[0,1]$ by $I$. Let $\equiv$ denote equivalence modulo $2$.

In what follows assume all manifolds are  oriented and equipped with basepoints; specific orientations and basepoints will usually be suppressed.

\subsection{Intersection numbers in 4-manifolds}\label{intersections}

The link homotopy invariants investigated in this paper are closely related to the algebraic ``intersection numbers'' $\lambda$ and $\mu$ introduced by Wall \cite{W}. For a more thorough exposition of the latter invariants in the four-dimensional setting, see Chapter 1 of \cite{FQ}, from which our definitions are based.

Suppose $A$ and $B$ are properly immersed, self-transverse 2-spheres or 2-disks in a connected 4-manifold $Y$. (By self-transverse we mean that self-intersections arise precisely as transverse double points.) Suppose further that $A$ and $B$ are transverse and that each is equipped with a path (a \emph{whisker}) connecting it to the basepoint of $Y$.
\mycomment{
    The geometric characterization of the invariants $\lambda$ and $\mu$, defined below, is that the intersection number $\lambda(A,B)$ vanishes if and only if the intersection points between $A$ and $B$ may be equipped with \emph{Whitney disks}, and the self-intersection number $\mu(A)$ vanishes if and only if the analogous statement holds for the self-intersection points of $A$ \cite[Section 1.7]{FQ}.
}
For an intersection point $x\in A\cap B$, let $\lambda(A,B)[x] \in \pi_1(Y)$ denote the homotopy class of a loop that runs from the basepoint of $Y$ to $A$ along its whisker, then along $A$ to $x$, then back to the basepoint of $Y$ along $B$ and its whisker. Define $\sign_{A,B}[x]$ to be $1$ or $-1$ depending on whether or not, respectively, the orientations of $A$ and $B$ induce the orientation of $Y$ at $x$. The Wall intersection number $\lambda(A,B)$ is defined by the sum
\[
    \lambda(A,B) = \mysum{x\in A\cap B}{} \sign_{A,B}[x]\lambda(A,B)[x]
\]
in the group ring $\Z[\pi_1(Y)]$, and is invariant under homotopy rel boundary of $A$ or $B$  \cite[Proposition 1.7]{FQ}, but depends on the choice of basepoint of $Y$ and the choices of whiskers and orientations. For an element $h=\sum_i n_i g_i$ in $\Z[\pi_1(Y)]$ ($n_i\in \Z$, $g_i\in \pi_1(Y)$),  define $\overline{h}\in \Z[\pi_1(Y)]$ by $\sum_i n_i \overline{g_i}$. From the definition it is readily verified that $\lambda(B,A) = \overline{\lambda(A,B)}$ and that the following observations, which we record for later reference, hold.

\begin{proposition}\label{prop:lambda-product}
If $x, y\in A\cap B$, then the product of $\pi_1(Y)$-elements $\lambda(A,B)[x]\overline{(\lambda(A,B)[y])}$ is represented by a loop that runs from the basepoint to $A$ along its whisker, along $A$ to $x$, then along $B$ to $y$, and back to the basepoint along $A$ and its whisker. Moreover, if $Y$ has abelian fundamental group, then this group element does not depend on the choice of whiskers and basepoint.\qed
\end{proposition}

\begin{proposition}\label{prop:lambda-restricted}
If $D_A\subset A$ is an  immersed 2-disk that is equipped with the same whisker  and oriented consistently with $A$, then for each $x\in D_A\cap B$ we have $\lambda(A,B)[x] = \lambda(D_A,B)[x]$ and $\sign_{A,B}[x]=\sign_{D_A,B}[x]$. \qed
\end{proposition}

The intersection numbers respect sums in the following sense. Suppose that $A$ and $B$ as above are 2-spheres, and suppose there is an embedded arc $\gamma$ from $A$ to $B$, with interior disjoint from  both. Let $\iota_A$ be a path that runs along the whisker for $A$, then along $A$ to the initial point of $\gamma$, and let $\iota_B$ be a path that runs from the endpoint of $\gamma$, along $B$ and its whisker to the basepoint of $Y$. Form the connect sum $A\# B$ of $A$ and $B$ along $\gamma$ in such a way that the orientations of each piece agree with the result. Equipped with the same whisker as $B$, the 2-sphere $A\# B$ represents the element\footnote{Since $A$ and $B$ are both whiskered, we permit ourselves to confuse them with their respective homotopy classes in $\pi_2(Y)$.} $A+gB$ in the $\Z[\pi_1(Y)]$-module $\pi_2(Y)$, where $g = [\iota_A\ast \gamma\ast \overline{\iota_B}]\in \pi_1(Y)$.  If $C$ is an immersed 2-disk or 2-sphere in $Y$ transverse to $A$ and $B$, then $\lambda(A+gB, C) = \lambda(A,C) + g\lambda(B, C)$. The additive inverse $-A\in \pi_2(Y)$ is represented by reversing the orientation of $A$.

Allowing $A$ again to be a self-transverse 2-disk or 2-sphere, the Wall self-intersection number $\mu(A)$ is defined as follows. Let $f_A:D\to Y$ be a map with image $A$, where $D=D^2$ or $S^2$.  
Let $x$ be a double point of $A$, and let $x_1$, $x_2$ denote its two preimage points in $D$. If $U_1$, $U_2$ are disjoint neighborhoods of $x_1$, $x_2$ in $\int D$, respectively, that do not contain any other double point preimages, then the embedded 2-disks $f_A(U_1)$ and $f_A(U_2)$ in $A$ are said to be two different \emph{branches} (or sheets) intersecting at $x$. Let $\mu(A)[x] \in \pi_1(Y)$ denote the homotopy class of a loop that runs from the basepoint of $Y$ to $A$ along its whisker, then along $A$ to $x$ through one branch $f_A(U_1)$, then along the other branch $f_A(U_2)$ and back to the basepoint of $Y$ along the whisker of $A$. (Such a loop is said to \emph{change branches} at $x$.) Define $\sign_{A}[x]$ to be $1$ or $-1$ depending on whether or not, respectively, the orientations of the two branches of $A$ intersecting at $x$  induce the orientation of $Y$ at $x$. In the group ring $\Z[\pi_1(Y)]$, let
\[
    \mu(A) = \mysum{x}{} \sign_{A}[x]\mu(A)[x],
\]
where the sum is over all such self-intersection points. (Note that it may sometimes be more convenient to write  $\mu(f_A)=\mu(A)$.)  For a fixed whisker of $A$, changing the order of the branches in the above definition replaces $\mu(A)[x]$ by its $\pi_1(Y)$-inverse, so $\mu(A)$ is only well-defined in the quotient $Q(Y)$ of $\Z[\pi_1(Y)]$, viewed as an abelian group, by the subgroup $\{a-\overline{a}:a\in \pi_1(Y)\}$. The equivalence class of $\mu(A)$ in this quotient group is  invariant under regular homotopy rel boundary of $A$. 
Note also that if the 4-manifold $Y$ has abelian fundamental group, then $\mu(A)$ does not depend on the choice of whisker.

Let $\self(A)\in \Z$ denote the signed self-intersection number of $A$. The \emph{reduced} Wall self-intersection number $\hat \mu(A)$ may be defined by
\[
    \hat \mu(A) =  \mu(A) - \self(A)\in Q(Y).
\]
It is an invariant of \emph{homotopy} rel boundary \cite[Proposition 1.7]{FQ}; this observation derives from the fact that non-regular homotopy takes the form of local ``cusp'' homotopies which may each change $\mu(A)$ by $\pm 1$ (see \cite[Section 1.6]{FQ}.) 

\mycomment{
    The geometric characterization of the invariants $\lambda$ and $\mu$ is that the intersection number $\lambda(A,B)$ vanishes if and only if the intersection points between $A$ and $B$ may be equipped with \emph{Whitney disks}, and the self-intersection number $\mu(A)$ vanishes if and only if the analogous statement holds for the self-intersection points of $A$ \cite[Section 1.7]{FQ}.
}

\subsection{The link homotopy invariants}\label{sec:linkhtpy}

We now recall the definitions of the link homotopy invariants $\sigma$ of Kirk \cite{Ki1} and $\omega$ of Li \cite{Li97}.

Let $f:S^2_+\cup S^2_-\to S^4$ be a link map. After a link homotopy (in the form of a perturbation) of $f$ we may assume the restriction $f_{\pm}=f|_{S^2_\pm}:S^2_{\pm}\to S^4\setminus f(S^2_{\mp})$ to each component is a  self-transverse immersion.  Let $X_-=S^4\setminus f(S^2_-)$
 and choose a generator $s$ for $H_1(X_-)\cong \Z$, which we write multiplicatively. 
For each double point $p$ of $f(S^2_+)$, let $\alpha_p$ be a simple circle on $f(S^2_+)$ that changes branches at $p$ and does not pass through any other double points. We call $\alpha_p$ an \emph{accessory circle} for $p$.  Letting $n(p) = \lk(\alpha_p, f(S^2_+))$,  one defines
\[
    \sigma_+(f) = \sum_{p} \sign(p) (s^{|n(p)|}-1)
\]
in the ring $\Z[s]$ of integer polynomials, where the sum is over all double points of $f(S^2_+)$, and to simplify notation we write $\sign(p) = \sign_{f_+(S^2)}[p]$.

Reversing the roles of $f_+$ and $f_-$, we similarly define $\sigma_-(f)$ and write $\sigma(f)=(\sigma_+(f), \sigma_-(f))\in \Z[s]\oplus \Z[s]$. Kirk showed in \cite{Ki1} that $\sigma$ is a link homotopy invariant,  and in \cite{Ki2} that if $f$ is link homotopic to a link map for which one component is embedded, then $\sigma(f)=(0,0)$.

Let $\rho:\pi_1(X_-)\to H_1(X_-)=\Z\langle s\rangle$ denote the Hurewicz map. Referring to the definition of $\mu$ in the preceding section as applied to the map $f_+:S^2_\pm\to X_-$,   observe that $\rho$ carries $Q(X_-)$ to the ring of integer polynomials $\Z[s]$ and Kirk's invariant $\sigma_+$ is given by
\begin{align}\label{eq:sigma-mu}
    \sigma_+(f) = \rho(\hat \mu(f_+)) \in \Z[s].
\end{align}

As in \cite{me2}, we say that $f$ is $\pm$-\textit{good} if $\pi_1(X_\mp)\cong \Z$ and the restricted map $f_\pm$ is a self-transverse immersion with $\self(f_\pm)=0$. We say that $f$ is \emph{good} if it is both $+$- and $-$-good. Equation \eqref{eq:sigma-mu} has the following consequence.
\begin{proposition}\label{prop:sigma-mu}
If $f$ is a $\pm$-good link map, then $\sigma_\pm(f)=\mu(f_\pm)$.
\end{proposition}

The invariant $\sigma_\pm(f)$  obstructs, up to link homotopy, pairing up double points of $f(S^2_\pm)$ with Whitney disks in $X_\mp$. While the essential purpose of Whitney disks is to embed (or separate) surfaces (see \cite[Section 1.4]{FQ}), our focus will be on their \emph{construction} for the purposes of defining certain invariants. In the setting of link maps, the  following standard result (phrased in the context of link maps) is the key geometric insight behind all the invariants we discuss in this paper and will find later application. 

\begin{lemma}\label{lem:linking-number-whitney-circle}
Let $f$ be a link map such that $f_-$ is a self-transverse immersion, and suppose $\{p^+,p^-\}$ are a pair of oppositely-signed double points of $f(S^2_-)$. Let $U$ and $V$ each be an embedded 2-disk neighborhood of $\{p^+,p^-\}$ on $f(S^2_-)$ such that $U$ and $V$   intersect precisely at these two points. On $f(S^2_-)$, let $\alpha^+, \alpha^-$ be loops based at $p^+, p^-$ (respectively) that leave along $U$ and return along $V$. Let $\gamma_U$, $\gamma_V$ be oriented paths in $U, V$ (respectively) that run from $p^+$ to $p^-$. Then the oriented loop $\gamma_U\cup \overline{\gamma_V}$ satisfies
\[
\pushQED{\qed}
    \bigl|\lk(\gamma_U\cup \overline{\gamma_V}, f(S^2_+))\bigr| = \bigl|\lk(\alpha^+, f(S^2_+)) - \lk(\alpha^-, f(S^2_+))\bigr|.\qedhere
\popQED
\]
\end{lemma}

 Wishing to obtain a ``secondary'' obstruction, in \cite{Li97} Li   proposed the following $(\Z_2\oplus \Z_2)$-valued invariant to measure intersections between $f(S^2_+)$ and the interiors of these disks.

Suppose $f$ is a $-$-good link map with $\sigma_-(f)=\mu(f_-)=0$.  The double points of $f(S_-^2)$  may then be labeled  $\{p_i^+,p_i^-\}_{i=1}^k$ so that $\sign(p_i^+)=-\sign(p_i^-)$ and $n_i:=|n(p_i^+)|=|n(p_i^-)|$; consequently, by Lemma \ref{lem:linking-number-whitney-circle} we may let $f^{-1}(p_i^+) = \{x_i^+, y_i^+\}$ and $f^{-1}(p_i^-) = \{x_i^-, y_i^-\}$ so that if $\gamma_i$ is an arc on $S^2_-$ connecting  $x_i^+$ to $x_i^-$ (and missing all other double point  preimages) and $\gamma_i'$ is an arc on $S^2_-$ connecting  $y_i^+$ to $y_i^-$ (and missing $\gamma_i$ and all other double point  preimages),  then the loop $f(\gamma_i)\cup f(\gamma_i')\subset f(S^2_-)$ is nullhomologous, hence nullhomotopic, in $X_+$. Let $U_i$ (respectively, $U_i'$) be a neighborhood of $\gamma_i$ (respectively, $\gamma_i'$) in $S^2_-$. The arcs $\{\gamma_i,\gamma_i'\}_{i=1}^k$ and neighborhoods $\{U_i, U_i'\}_{i=1}^k$ may be chosen so that the collection $\{U_i\}_{i=1}^k\cup \{U_i'\}_{i=1}^k$ is mutually disjoint,  and so that the resulting \emph{Whitney circles} $\{f(\gamma_i\cup \gamma_i')\}_{i=1}^k$ are mutually disjoint, simple circles in $X_+$ such that each bounds an immersed \emph{Whitney disk} $W_i$ in $X_+$ whose interior is transverse to $f(S^2_-)$. 

Since  the two branches $f(U_i)$ and $f(U_i')$ of $f(S^2_-)$ meet transversely at $\{p_i^+, p_i^-\}$, there are a pair of smooth vector fields $v_1, v_2$ on $\d W_i$ such that $v_1$ is tangent to $f(S^2_-)$ along $\gamma_i$ and normal to $f(S^2_-)$ along $\gamma_i'$, while $v_2$ is normal to $f(S^2_-)$ along $\gamma_i$ and tangent to $f(S^2_-)$ along $\gamma_i'$. Such a pair defines a normal framing of $W_i$ on the boundary. We say that $\{v_1,v_2\}$ is a \emph{correct framing} of $W_i$, and that $W_i$ is \emph{framed}, if the pair extends to a normal framing of $W_i$.  By boundary twisting $W_i$ (see page 5 of \cite{FQ}) if necessary, at the cost of introducing more interior intersection points with $f(S^2_-)$, we can choose the collection of Whitney disks $\{W_i\}_{i=1}^k$ such that each is correctly framed.

Let $1\leq i\leq k$. To each point of intersection $x  \in  f(S^2_-)\,\cap \,\int W_i$, let $\beta_x$ be a loop that first goes along $f(S^2_-)$ from its basepoint to $x$, then along $W_i$ to $f(\gamma_i')\subset \d W_i$, then back along $f(S^2_-)$ to the  basepoint of $f(S^2_-)$ and let $m_i(x)=\lk(f(S^2_+), \beta_x)$. Let
\begin{align}\label{eq:Jix}
    {\cL}^-_i(x)= n_i + n_im_i(x) + m_i(x) \mod 2; 
\end{align}
summing over all such points of intersection, let
\[
    {\cL}^-(W_i) = \mysum{x\, \in\,  f(S^2_-)\,\cap \,\int W_i}{} {\cL}^-_i(x) \mod 2.
\]
Then Li's $\omega_-$-invariant applied to $f$  is defined by
\begin{align}\label{eq:omega}
    \omega_-(f) = \mysum{i=1}{k}\; \cL^-(W_i) \mod 2.
\end{align}

We record  two  observations about this definition for later use. The latter is a special case of Proposition  \ref{prop:lambda-product}.

\begin{remark}\label{rem:omega-compute}
If $n_i$ is odd then $\cL^-(W_i) = f(S^2_-)\cdot \int W_i \mod 2$, while if $n_i$ is even then
\[
    \cL^-(W_i) = \mysum{x\, \in\,  f(S^2_-)\,\cap \,\int W_i}{} m_i(x).
\]
\end{remark}

\begin{remark}\label{rem:omega-product}
Suppose $x, y \in f(S^2_-)\cap \int W_i$, and let $\beta$ be a loop that runs from $x$ to $y$ along $f(S^2_-)$, then back to $x$ along $\int W_i$. We have
\[
    m_i(x) + m_i(y) = \lk(f(S^2_+), \beta) \textup{ mod } 2.
\]
\end{remark}


Now suppose $f$ is an arbitrary link map with $\sigma_-(f)=0$. By standard arguments 
we may choose a $-$-good link homotopy representative $f'$ of $f$, and  $\omega(f)$ is defined by setting $\omega_-(f)=\omega_-(f')$. By a result in \cite{me2}, in \cite{Li97} it was shown that this defines an invariant of \emph{link homotopy}; Theorem \ref{thm:mainresult} gives a new proof. By  interchanging the roles of $f_+$ and $f_-$ (and instead assuming $\sigma_+(f)=0$), we obtain $\omega_+(f)\in \Z_2$ similarly, and write $\omega(f) = (\omega_+(f), \omega_-(f))$.

Based on similar geometric principles, Teichner and Schneiderman \cite{ST} defined a secondary obstruction with respect to the homotopy invariant $\mu$. When adapted to the context of link homotopy, however, their invariant reduces to $\omega$ \cite{me2}.

 \subsection{Surgering tori to 2-spheres}\label{surger}


A common situation to arise in this paper is the following. Suppose we have a  torus (punctured torus, respectively) $T$ in the 4-manifold $Y$,  on which we wish to perform surgery along a curve so as to turn it into a 2-sphere  (2-disk, respectively) whose  self-intersection and intersection numbers may be calculated.  Our device for doing so is the following lemma, which is similar to \cite[Lemma 4.1]{me2} and so its proof is omitted. A similar construction may be found at the bottom of page 86 of \cite{FK}.

If $D_A\subset A$ is an immersed 2-disk, let $\mu(A)|_{D_A}$ denote the contribution to $\mu(A)$ due to self-intersection points on $D_A$.

\begin{lemma}\label{lem:surgery}
Suppose that $Y$ is a codimension-$0$ submanifold of $S^4$ and $\pi_1(Y)$ is abelian. Let $B$ be a properly immersed 2-disk or 2-sphere in $Y$, and suppose  $T$ is an embedded torus {(}or punctured torus{)} in $Y\setminus \int B$. Let  $\delta_0$ be a simple, non-separating curve on $T$, let $\delta_1$ be a normal pushoff of $\delta_0$ on $T$ and let $\hat T$ denote the annulus on $T$ bounded by $\delta_0\cup \delta_1$. Suppose there is a map
\[
    J: D^2\times I\to   Y
\]
such that $J(D^2\times I)\cap T = J(\d D^1\times I) = \hat T$ and $J(\d D^1\times \{i\}) = \delta_i$ for $i=0,1$. Then, after a small perturbation,
\[
    S = (T \setminus \int \hat T) \mycup{\delta_0\cup \delta_1}{} J(D^2\times \{0,1\})
\]
is a properly immersed, self-transverse 2-sphere (or 2-disk, respectively) in $Y$ such that
\begin{enumerate}[label=\textup{(\roman*)}, ref=\textup{(\roman*)},align=CenterWithParen]
\item $\lambda(S,B) = (1-[\delta_2])\lambda(D,B)$ in $\Z[\pi_1(Y)]$, and
\item $\mu(S) \equiv \mu(S)|_{D} + \mu(S)|_{D'} +(D\cdot D') [\delta_2]\mycomment{\overline{[\delta_2]}}$ in $Q(Y)\hskip-0.1cm\mod 2$,
\end{enumerate}
where\footnote{Note that we shall sometimes exclude curly brackets for a one-point set that occurs in a Cartesian product.} $D = J(D^2\times 0)\subset S$ 
 is oriented consistently with and shares  a whisker with $S$, $D'=J(D^2\times 1)$, and $\delta_2$ is a dual curve to each of $\delta_0$ and $\delta_1$ on $T$ such that $\delta_2\cap \hat T$ is a simple arc running from $\delta_0$ to $\delta_1$.\qedhere
\end{lemma}

The hypotheses of this lemma will frequently be encountered in the following form. We have $T$ and the nullhomotopic curve $\delta_0$; we then choose an immersed 2-disk $D\subset Y$ bounded by $\delta_0$ and let $J$ be its ``thickening'' along a section (which is not necessarily non-vanishing) obtained by extending over the 2-disk a normal section to $\delta_0$ that is tangential to $T$.

\section{Proof of Theorem \ref{thm:mainresult}}\label{sec:proof}

Let us first outline the steps of our proof. Up to link homotopy, one component $f(S^2_-)$ of a link map $f$ is unknotted, immersed, and in Section \ref{sec:unknotted} we exploit this to construct a collection of mutually disjoint, embedded, framed Whitney disks $\{V_i\}_i$ for $f(S^2_-)$ with interiors in $S^4\setminus f(S^2_-)$, such that each has nullhomotopic boundary in the complement of $f(S^2_+)$. We show how these disks (along with disks  bounded by accessory circles for $f_-$) can be used to construct  2-sphere generators of  $\pi_2(X_-)$. The algebraic intersections between $f(S^2_+)$ and these 2-spheres are then computed in terms of the intersections between $f(S^2_+)$ and the aforementioned  disks.

In Section \ref{sec:omega} we surger the disks $\{V_i\}_i$ so to exchange their intersections with $f(S^2_+)$ for intersections with $f(S^2_-)$. In this way we obtain \emph{immersed}, framed Whitney disks $\{W_i\}_i$ for $f_-$ in $S^4\setminus f(S^2_+)$, such that the algebraic intersections between $f(S^2_-)$ and $W_i$, measured by  $\omega_-(f)$, are related to the algebraic intersections between $f(S^2_+)$ and $V_i$.

In Section \ref{sec:main3} we complete the proof by combining the results of these two sections: the intersections  between $f(S^2_+)$ and  generators of $\pi_2(X_-)$ of the former section are related to $\sigma_+(f)$, which by the latter section can be related to $\omega_-(f)$.

\subsection{Unknotted immersions and Whitney disks in $X_-$}\label{sec:unknotted}


A notion of unknottedness for surfaces in 4-space was introduced  by Hosokawa and Kawauchi  in \cite{Hos}. 
A connected, closed, orientable surface in $\R^4$ is said to be unknotted if it bounds an embedded 3-manifold in $\R^4$ obtained by attaching (3-dimensional)  $1$-handles to a $3$-ball. They showed that by attaching (2-dimensional) 1-handles only, any embedded surface in $\R^4$ can be made unknotted. 
Kamada \cite{K} extended their definition to \emph{immersed} surfaces in $\R^4$ and gave a notion of equivalence for such immersions.   In that paper it was shown that an immersed 2-sphere in $\R^4$ can be made equivalent, in this sense, to an unknotted, immersed one by performing (only) \emph{finger moves}.

It was noted in \cite{me2} that we may perform a link homotopy to ``unknot'' one immersed component of a link map (see Lemma \ref{lem:good}).  The algebraic topology of the complement of this unknotted (immersed) component is greatly simplified,  making the computation of the invariants defined in the previous section more tractable.

\begin{figure}%
    \centering
    \subfloat[The immersed 2-disk $D^+$ with one double point $r^+$.]{%
    \label{fig:fig-cusp-positive}%
    \includegraphics[width=\linewidth]{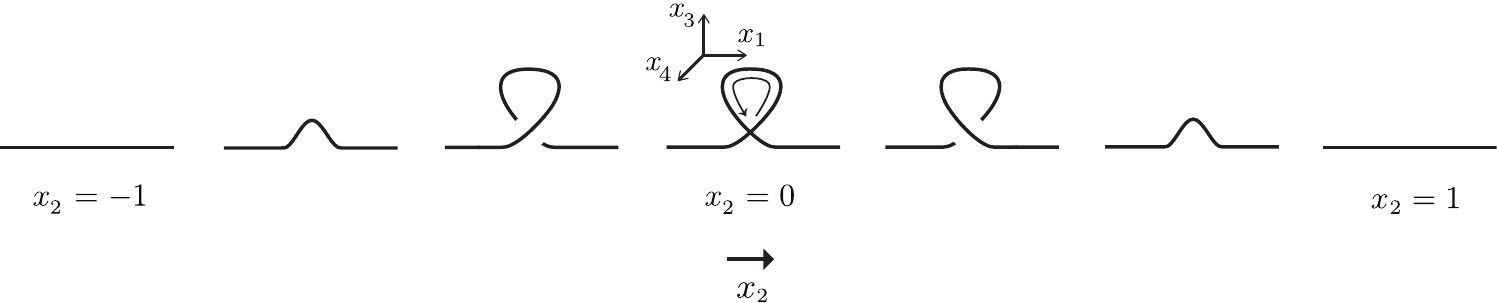}%
    }%
    \vskip 1cm
    \subfloat[The immersed 2-disk $D^-$ with one double point $r^-$.]{%
    \label{fig:fig-cusp-negative}%
    \includegraphics[width=\linewidth]{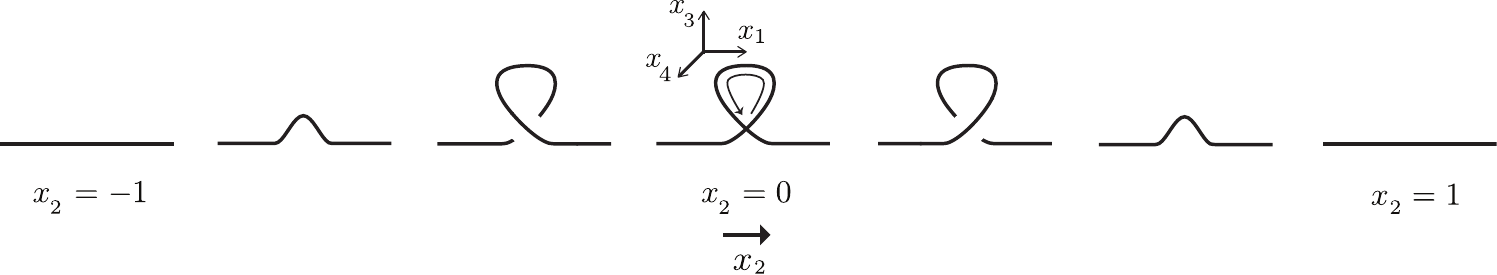}%
    }%
    \caption{}%
    \label{fig:fig-cusps}%
\end{figure}

Let us begin with a precise definition of an unknotted, immersed 2-sphere. To do so we construct cusp regions which have certain symmetry properties;  our justification for these specifications is to follow.
Using the moving picture method, Figures \ref{fig:fig-cusps}(a), (b) (respectively) illustrate properly immersed, oriented 2-disks $D^+$, $D^-$ (respectively) in $D^4$, each with precisely one double point $r^+$, $r^-$  (respectively) of opposite sign. In those figures we have  indicated coordinates $(x_1,x_2,x_3,x_4)$ of $D^4$; our choice of the $x_2$-ordinate to represent ``time'' is a compromise between  ease of illustration and ease of subsequent notation.

As suggested by these figures, we construct $D^\pm$ so that it has boundary $\d D^\pm = \d D^2\times 0\times 0$ and so that it intersects $D^1\times 0\times D^1\times D^1$ in an arc lying in the plane $D^1\times0\times D^1\times 0$. Further, letting $\theta^\pm$ denote the loop on $D^\pm$ in this plane that is  based at $r^\pm$ and oriented as indicated in those figures, we have that the reverse loop $\overline{\theta^\pm}$ is given by\footnote{More precisely, we can choose a parameterization $\theta^\pm:I\to D^1\times 0\times D^1\times 0$ and paths $\theta_i^\pm:I\to D^1$, $i=1,3$, so that for $t\in I$ we have $\theta^\pm(t)=(\theta_1^\pm(t), 0,\theta_3^\pm(t), 0)$ and hence
  \[
    \Sigma(\theta^\pm(t)) =  (-\theta_1^\pm(t), 0,\theta_3^\pm(t), 0) = \overline{\theta^\pm}(t).
\]}
\begin{align}
        \Sigma \circ \theta^\pm &= \overline{\theta^\pm},\label{eq:propertiesofSigma-1}
\end{align}
where $\Sigma$ is the orientation-preserving self-diffeomorphism of $D^4$ given by $(x_1,x_2,x_3,x_4)\mapsto (-x_1, -x_2, x_3, x_4)$.
Lastly, we may suppose that
\begin{align}
\Sigma(D^\pm) &= D^\pm.\label{eq:propertiesofSigma-2}
\end{align}

After orienting $D^4$ appropriately, the immersed 2-disk $D^+$ ($D^-$, respectively) has a single, positively (negatively
, respectively) signed double point $r^+$ ($r^-$, respectively), and $\theta^\pm$ is an oriented loop on $D^\pm$ based at $r^\pm$ which changes branches there. (Note also that in this construction we may suppose that $D^-$ is the image of $D^+$ under the orientation-reversing self-diffeomorphism of $D^4$ given by $(x_1,x_2,x_3,x_4)\mapsto (-x_1, x_2, x_3, x_4)$.) We call $D^+$ and $D^-$ \emph{cusps}.

Roughly speaking, an unknotted immersion is obtained from an unknotted embedding by ``grafting on'' cusps of this form. The purpose of the above specifications is so that, by a manoeuvre resembling  the Disk Theorem, we may more conveniently move these cusps around on the 2-sphere so that accessory circles of the form $\theta^\pm$ are permuted and perhaps reversed as oriented loops.

Formally, let $d\geq 0$,  let $\eps:\{1,2,\ldots, d\}\to \{+,-\}$ be a map which associates a $+$ sign or a $-$ sign to each $i\in \{1,2,\ldots, d\}$, and write $\eps_i=\eps(i)$. Let $\U$ denote the image of an  oriented, unknotted embedding $S^2\to S^4$; that is, an embedding that extends to the 3-ball (which is unique up to ambient isotopy). Suppose there are a collection of mutually disjoint, equi-oriented embeddings $b_i:D^4\to S^4$, $i=1,\ldots, d$, such that
\[
    b_i^{-1}(\U) = D^2\times 0\times 0
\]
for each $1\leq i\leq d$. By removing the interiors of the 2-disks $\{b_i(D^2\times 0\times 0)\}_i$ from $\U$ and attaching, for each $1\leq i\leq d$, the cusp $b_i(D^{\eps_i})$ along $b_i(\d D^2\times 0\times 0)$, we obtain an \emph{unknotted, immersed} 2-sphere in $S^4$:
\begin{align}\label{defn:unknotted}
    \U_\eps = [\U\setminus \mycup{i=1}{d} \int  b_i(D^4)] \, \cup\,\mycup{i=1}{d} b_i(D^{\eps_i}).
\end{align}

Note that we use the function $\eps$ in \eqref{defn:unknotted} only for convenience of notation in the proofs that follow. Since the embeddings $\{b_i\}_i$ can always be relabeled, one sees that the definition of $\U_\eps$ depends precisely on the choice of unknotted 2-sphere $\U$, the embeddings $\{b_i\}_i$, and two non-negative integers $d_+$ and $d_-$, where $d_\pm$ is the number of $1\leq i\leq d$ such that $\eps_i=\pm$. 

The following lemma will allow us to perform an ambient isotopy of $S^4$ which carries the model $\U_\eps$ back to itself such that accessory circles are permuted (and perhaps reversed in orientation) in a prescribed manner.

\begin{lemma}\label{prop:exchange-all-b-i}
Let $\rho$ be a permutation of $\{1,2,\ldots, d\}$. For each $1\leq i\leq d$, let $\mu_i\in \{-1,1\}$. There is an ambient isotopy $\hat\varphi:S^4\times I\to S^4$  such that $\hat\varphi_1$ fixes
\[
    \U_\eps\setminus \mycup{i=1}{d} \int  b_i(D^4)
\]
set-wise\mycomment{\footnote{A map $f$ is said to fix a set $A$ set-wise if $f(A)=A$.}}  and, for each $1\leq i\leq d$,
\[
    \hat\varphi_1\circ b_i(x,y) = b_{\rho(i)}(\mu_i \, x, y)
\]
for all $(x,y)\in D^2\times D^2$.  
\mycomment{
    In particular,  $\hat\varphi_1$ sends $b_j(D^{\eps_j})$ to $b_{\rho(j)}(D^{\eps_j})$, so if $\rho$ preserves parity then $\hat\varphi_1(e_d)=e_d$.
}%
In particular, if $\eps_{\rho(i)}=\eps_i$, then $\varphi_1(\U_\eps\cap b_i(D^4)) = \U_\eps\cap b_{\rho(i)}(D^4)$ and
\begin{align*}
\hat \varphi_1\circ b_i\circ \theta^{\eps_i} = \begin{cases} b_{\rho(i)}\circ \theta^{\eps_i} & \text{ if $\mu_i=1$,}\\[0.2cm] b_{\rho(i)}\circ \overline{\theta^{\eps_i}} & \text{ if $\mu_i=-1$.} \end{cases}
\end{align*}
\end{lemma}

The proof consists of using the Disk Theorem \cite[Corollary 3.3.7]{Kosinski} to transport 4-ball neighborhoods of the cusps around the 2-sphere, and is deferred to Appendix \ref{app:prop:FM-constructVis-withu}. We proceed instead to apply the lemma to equip one component of a link map, viewed as an immersion into the 4-sphere, with a particularly convenient collection of mutually disjoint, embedded, framed Whitney disks.


For this purpose it will be useful to give a particular construction of an unknotted immersion of a 2-sphere in $S^4$ with $d\geq 0$ pairs of opposite-signed double points.

\subsubsection{A model, unknotted immersion}

\begin{figure}[h]
    \centering
        \includegraphics[width=.25\linewidth]{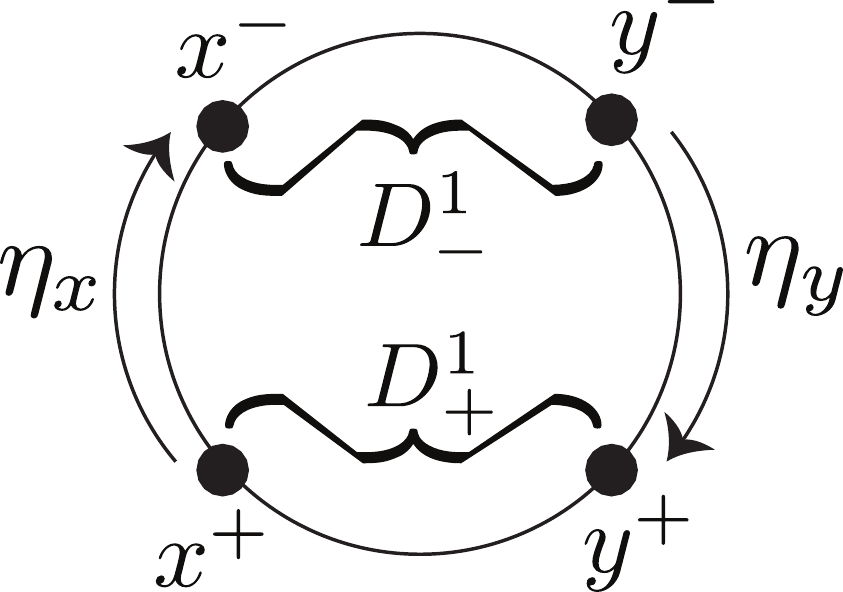}
        \caption{}
        \label{fig:circle}
    \end{figure}

For $A\subset \R^3$ and real numbers $a<b$, write $A[a,b]=A\times [a,b]\subset \R^3\times \R$, and $A[a] = A\times a$. Choose an increasing sequence
\[
        0 = t_1^-<t_1^+<t_2^-<t_2^+<\ldots<t_{d-1}^+<t_d^-<t_d^+=1,
\]
and for each $1\leq i\leq d$, write $I_i = [t_i^-,t_i^+]$ and let $t_i=(t_i^+ + t_i^-)/2$. On the unit circle, oriented clockwise, let $x^+,x^-,y^-$ and $y^+$ be distinct, consecutive points and let $D^1_+$, $D^1_-$ be disjoint neighborhoods of $\{x^+,y^+\}$, $\{x^-,y^-\}$, respectively. Let $ \hat \alpha^\pm:I\to D^1_\pm$ be a path in $D^1_\pm$ running from $x^\pm$ to $y^\pm$; let $\eta_x$ and $\eta_y$ be  simple paths on $S^1$ running $x^+$ to $x^-$ and from $y^-$ to $y^+$, respectively. Let $\hat \eta_x$ and $\hat \eta_y$ be disjoint neighborhoods of $\eta_x$ and $\eta_y$ in $S^1$, respectively.
    See Figure \ref{fig:circle}. 

\begin{figure}[h]
\centering
    \includegraphics[width=\linewidth]{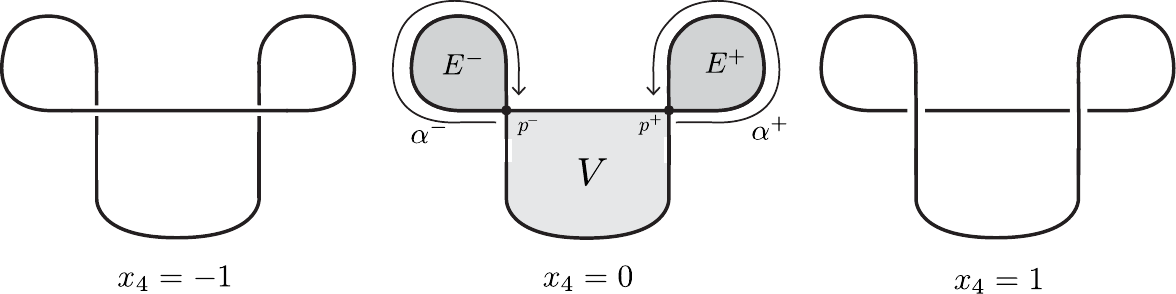}
    \caption{}
    \label{fig:fig-standard-annulus-E-alpha-p-V}
\end{figure}

For each $1\leq i\leq d$, let $\hat \Theta_i:D^1\to I_i$ be a linear map such that $\hat \Theta_i(0)=t_i$ and $\hat \Theta_i(\pm 1)=t_i^\pm$, and let $\Theta_i:D^3\times D^1\to D^3\times I_i$ be the map $\Theta_i(x,t)=(x,\hat \Theta_i(t))$.  Let $G:S^1\times D^1\to D^3\times D^1$ be an oriented, self-transverse immersion with image as shown in Figure \ref{fig:fig-standard-annulus-E-alpha-p-V} (ignoring the shadings), with two double-points $p^\pm=G(x^\pm,0)=G(y^\pm, 0)$, such that $G(x,t)\subset D^3[t]$ for each $x\in S^1$ and $t\in D^1$. Then  $\alpha^\pm = G(\hat\alpha^\pm\times 0)$ is an oriented loop on $G(S^1\times D^1)$ based at $p^\pm$ which changes branches there. Note that $G$ is the trace of a regular homotopy from the circle in $D^3$ to itself that Figure \ref{fig:fig-standard-annulus-E-alpha-p-V} illustrates.

\mycomment{Now, $G(S^1\times I)$ is obtained by gluing the two opposite-signed cusps $G(D^1_+\times I)$ and $G(D^1_-\times I)$ along $G(\d D^1_+\times I)$. }

For each $1\leq i\leq d$, define a map $G_i:S^1\times I_i\to D^3\times I_i$ by
\[
    G_i(\Theta_i(x,t)) = \Theta_i(G(x,t))
\]
for $(x,t)\in S^1\times D^1$. \mycomment{note that $S^1\times [t_i^-,t_i^+] = \Theta_i(S^1\times D^1)$.}
Write the 2-sphere as the capped off cylinder
 \begin{align}
    S^2 = \bigl(D^2\times\{-1,1\}\bigr) \cup \mycup{i=1}{d} \Theta_i(S^1\times I) 
 \end{align}
 in $D^3\times D^1$, and define a map $\hat u_d:S^2\to D^3\times D^1\subset S^4$ by the identity on $D^2\times \{-1,1\}$ and by  $G_i$ on $S^1\times I_i=\Theta_i(S^1\times D^1)$.

\mycomment{By the above construction there are equi-oriented diffeomorphisms $b_i^\pm:D^4\to \Theta_i(B^3_\pm\times I) = B^3_\pm\times[t_i^-,t_i^+]$ given by $b_i^\pm = \Theta_i\circ \Gamma$ taking $D^\pm$ to $\Theta_i(G(D^1_\pm\times I))$\mycomment{$=G_i(\Theta_i(D^1_\pm\times I))}, so by Proposition \ref{prop:permute-cusps}, we have the following.}


After smoothing corners, $\hat \U_d = \hat u_d(S^2)$ is an immersed 2-sphere in $S^4$. Let $\alpha_i^\pm$ be the oriented loop on $\hat \U_d$ given by $\alpha_i^\pm(s) =  \Theta_i(G(\hat \alpha^\pm(s),0)) = \Theta_i(\alpha^\pm(s))$ for $s\in I$. Observe that $\alpha_i^\pm$ is based at the $\pm$-signed double point $p_i^\pm = \Theta_i(G(x^\pm,0))=\Theta_i(G(y^\pm,0))$ and changes branches there.

\begin{figure}[h]
\centering
    \includegraphics[width=0.7\linewidth]{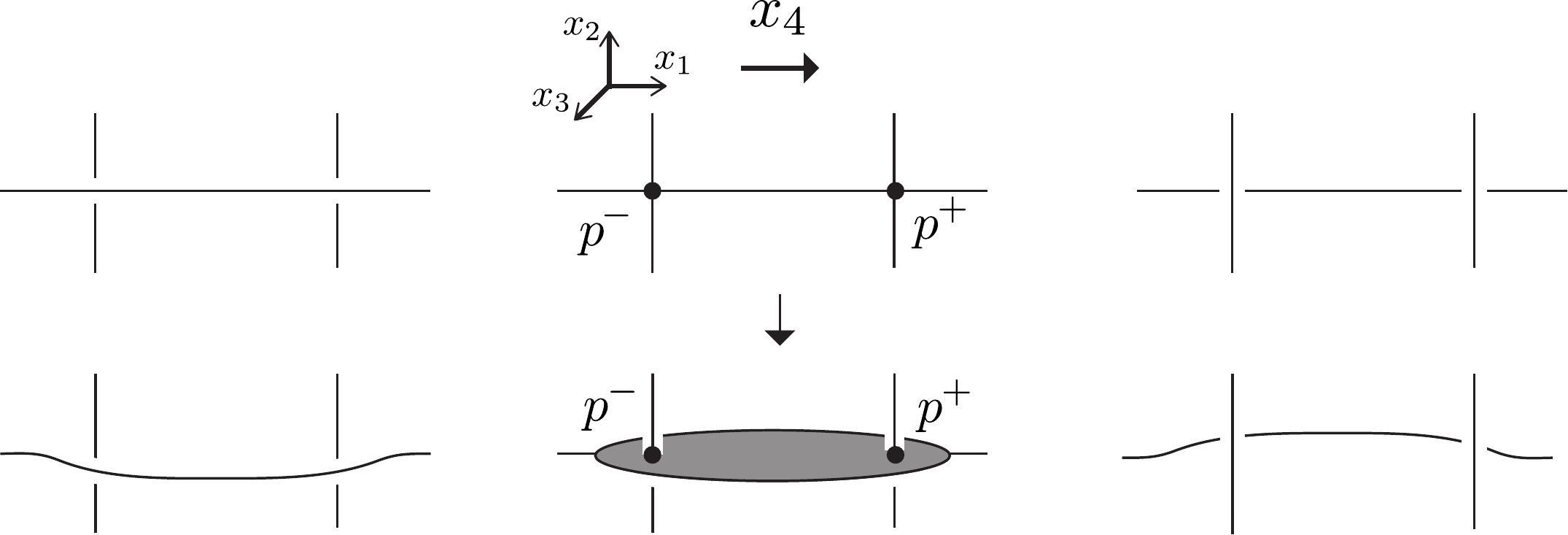}
    \caption{}
    \label{fig:fig-fatten}
\end{figure}

Referring to Figure \ref{fig:fig-standard-annulus-E-alpha-p-V}, let $V[0]\subset D^3[0]$ ($V\subset D^3$) be the obvious, embedded Whitney disk for the immersed annulus $G(S^1\times D^1)$ in $D^3\times D^1$, bounded by $G((\eta_x\cup \eta_y)\times 0)$. For arbitrarily small $\eps>0$, by pushing a neighborhood of $G(\eta_y\times 0)$ in $G(\hat\eta_y\times (-\eps,\eps))$ into $D^3\times 0$ (as in Figure \ref{fig:fig-fatten}), we may assume that the Whitney disk is framed: the constant vector field $(0,0,1,0)$ that points out of the page in each hyperplane $\R^3[t_0]$ and the constant vector field $(0,0,0,1)$ are a correct framing.  Thus  the maps $\{\Theta_i\}_{i=1}^d$ carry $V[0]$ to a complete collection of mutually disjoint, embedded, framed Whitney disks $V_i\subset \R^3[t_i]$, $i=1,\ldots, d$, for $\hat \U_d$ in $S^4$. In particular, the boundary of $V_i$ (equipped with an orientation) is given by
\begin{align}\label{eq:bdy-V-i}
    \d V_i = \Theta_i(G((\eta_x\cup \eta_y)\times 0)) = G_i(\eta_x\times 0) \cup G_i(\eta_y\times 0).
\end{align}

By Lemma 4.2 of \cite{me2}, we may assume after a link homotopy that one component of a link map is of the form of this model of an unknotted immersion.

\begin{lemma}\label{lem:good}
A link map $f$ is link homotopic to a good link map $g$ such that $g(S^2_-)=\hat \U_d$ for some non-negative integer $d$. \qed
\end{lemma}

We proceed to generalize the proof of  \cite[Lemma 4.4]{me2} to construct representatives of a basis of $\pi_2(X_+)$ and compute their algebraic intersections with $f(S^2_+)$ in terms of the algebraic intersections between $f(S^2_+)$ and $\{V_i\}_i$.  While parts (i) and (iii) of Proposition \ref{lem:2-spheres} may be deduced directly from that paper, we include the complete proof for clarity.

\begin{proposition}\label{lem:2-spheres}
Let $f$ be a good link map such that $f(S^2_-)=\hat \U_d$. Equip $f(S^2_+)$ with a whisker in $X_-$ and fix an identification of $\pi_1(X_-)$ with $\Z\langle s\rangle$ so as to write $\Z[\pi_1(X_-)] = \Z[s,s^{-1}]$.  Then $\pi_2(X_-)\cong (\hskip -0.03cm\underset{i=1}{\overset{d}{\oplus}}\Z)[s,s^{-1}]\mycomment{(\Z[\Z])^d}$ and there is a $\Z[s,s^{-1}]$-basis represented by mutually disjoint, self-transverse, immersed, whiskered 2-spheres $\{A_i^+, A_i^-\}_{i=1}^d$ in $X_-$ with the following properties. For each $1\leq i\leq d$, there is an integer Laurent polynomial $q_i(s)\in \Z[s,s^{-1}]$ such that
\begin{enumerate}[label=\textup{(\roman*)}, ref=\textup{(\roman*)},align=CenterWithParen]
    \item $n_i:= q_i(1) = \lk(f(S^2_+),\alpha_i^+)$,
    \item $\lambda(A_i^\pm, A_i^{\pm\prime}) = s+s^{-1} \mod 2$, 
    \item $\lambda(f(S^2_+), A_i^+) = (1-s)^2q_i(s)$, and
    \item $\lambda(f(S^2_+), A_i^-) - \lambda(f(S^2_+), A_i^+) = (1-s)^2\lambda(f(S^2_+),V_i^c)$,
\end{enumerate}
where $V_i^c$ is a 2-disk in $X_-$ obtained from $V_i$ by removing a collar in $X_+$, and $A_i^{\pm\prime}$ is the image of a section of the normal bundle of $A_i^\pm$. Moreover, if for any $1\leq j\leq d$ and $\eps\in \{-, +\}$ the loop $\alpha_j^\eps$ bounds a 2-disk in $S^4$ that intersects $f(S^2_+)$ exactly once, then we may choose $A_j^{\eps}$ so that
\[
    \lambda(f(S^2_+), A_j^{\eps}) = (1-s)^2.
\] \mycomment{ I might want to use this to show Tau=Omega.}
\end{proposition}

\begin{proof} 
For each $1\leq i\leq d$, the double points $p_i^+$ and $p_i^-$ of $f(S^2_-)$ lie in $D^3[t_i]$ and  $f(S^2_-)$ intersects the 4-ball $D^3[t_i^-, t_i^+]\subset D^3\times D^1$ precisely along $\Theta_i(G(S^1\times D^1))$. 
In what follows we shall denote $\hat G = G(S^1\times D^1)$; note that for $t\in D^1$, $\hat G\cap (D^3[t])=G(S^1\times t)$ and $\Theta_i(G(S^1\times t)) = f(S^2_-)\cap (D^3[\hat \Theta_i(t)])$.  Observe that by the construction of the annulus $G$ we may assume there are integers  $-1<a^-<b^-<0<a^+<b^+<1$ such that $G(S^1\times [a^\pm,b^\pm])=G_0[a^\pm, b^\pm]$ for some (embedded circle) $G_0\subset D^3$.

Let $E^\pm$ denote the 2-disk bounded by $\alpha^\pm$ in $D^3[0]$, and choose a 3-ball $N^\pm$ so that $N^\pm[-1,1]$ is a 4-ball neighborhood of $p^\pm$ and $\Theta_i(N^\pm[-1,1])$  is disjoint from $f(S^2_+)$. There is an embedded torus $T^\pm$ in $N^\pm[-1,1]\setminus \hat G$ that intersects $E^\pm$ exactly once; see Figure \ref{fig:fig-standard-annulus-E-T-intersection} for an illustration of  $T^+$ and $T^-$ in $D^3\times D^1$. The torus $T^\pm$ appears as a cylinder in each of 3-balls $N^\pm[-1]$ and $N^\pm[1]$, and as a pair of circles in $N^\pm[t]$ for $t\in (-1,1)$.
\newlength{\myheight}%
\setlength{\myheight}{0.25\textwidth}%
\begin{figure}[h]
\centering
    \includegraphics[width=\linewidth]{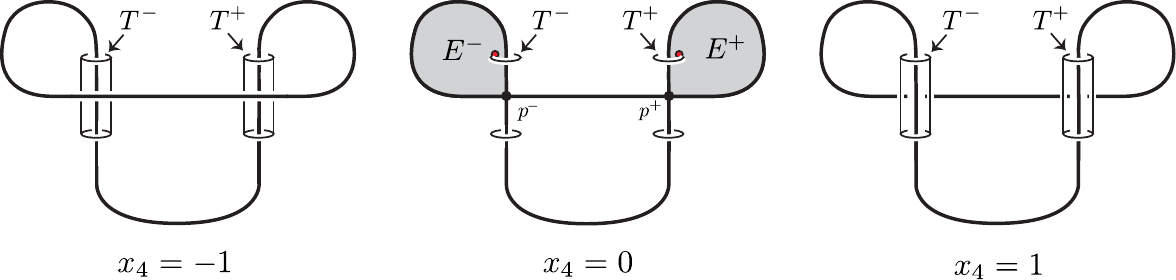}
    \caption{}
    \label{fig:fig-standard-annulus-E-T-intersection}
\end{figure}
For each $1\leq i\leq d$, let $T_i^\pm = \Theta_i(T^\pm)\subset N^\pm[t_i^-,t_i^+]\setminus f(S^2_-)$.  By Alexander duality the linking pairing
\[
    H_2(X_-)\times H_1(f(S^2_-))\to \Z
\]
defined by $(R,\upsilon)\mapsto R\cdot\Upsilon$, where $\upsilon=\d \Upsilon\subset S^4$, is nondegenerate. Thus, as the loops $\{\alpha_i^+,\alpha_i^-\}_{i=1}^d$ represent a basis for $H_1(f(S^2_-))\cong \Z^{2d}$, we have that $H_2(X_-) \cong \Z^{2d}$ and (after orienting) the so-called \textit{linking tori}  $\{T_i^+,T_i^-\}_{i=1}^d$ represent a basis.  
We proceed to apply the construction of Section \ref{surger} (twice, successively) to turn these tori into 2-spheres.

Let $\widehat \Delta^\pm\subset N^\pm$ be the embedded 2-disk  so that $\widehat \Delta^\pm[a^\pm]$ appears in $N^\pm[a^\pm]$ as in Figure \ref{fig:fig-delta-hat-Delta-gamma-plus-minus-gamma-tub-nbd-no-as}, and let $\hat \delta^\pm =  \d\widehat \Delta^\pm$. The disk $\widehat \Delta^\pm[a^\pm]$  intersects $\hat G$ at two points, which are the endpoints of an arc of the form $\gamma^\pm[a^\pm]$  ($\gamma^\pm\subset G_0$) in $G(S^1\times {a^\pm})$, also illustrated. 

In  Figure \ref{fig:fig-delta-hat-Delta-gamma-plus-minus-gamma-tub-nbd-no-as} we have also illustrated in $D^3[a^\pm]$  the restriction of a  tubular neighborhood   of $\hat G$  to $\gamma^\pm[a^\pm]$, which we may write in the form $(\gamma^\pm\times D^2)[a^\pm]$ and choose so that $\Theta_i$ carries $(\gamma^\pm\times D^2)[a^\pm]$ to a tubular neighborhood of $f(S^2_-)$ restricted to $\Theta_i(\gamma^\pm[a^\pm])$  in $S^4\setminus f(S^2_+)$. 
\begin{figure}[h]
\centering
    \includegraphics[width=.9\linewidth]{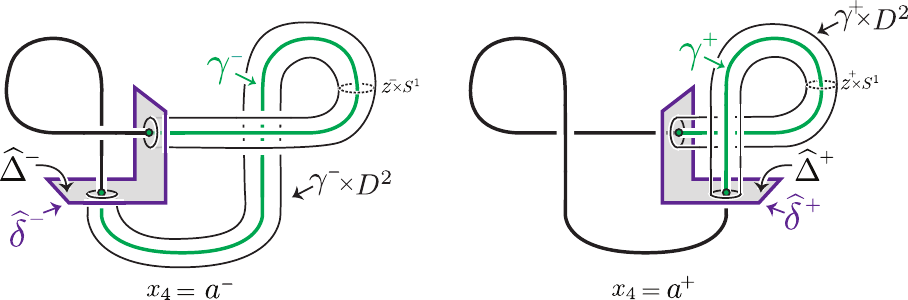}
    \caption{}
    \label{fig:fig-delta-hat-Delta-gamma-plus-minus-gamma-tub-nbd-no-as}
\end{figure}
Let $\Delta^\pm$ be the embedded punctured torus in $D^3$ obtained from $\hat \Delta^\pm$ by attaching a 1-handle along $\gamma^\pm$; that is, let
\[
    \Delta^\pm = \Bigl[\hat\Delta^\pm \setminus (\d \gamma\times \int D^2)\Bigr] \mycup{\d\gamma^\pm\times S^1}{} (\gamma^\pm\times S^1).
\]
Note that $\Delta^\pm$ has boundary $\hat \delta^\pm$.

Let $\beta^\pm$ be a loop on $\Delta^\pm$ formed by connecting the endpoints of a path on $\gamma^\pm\times S^1$ by an arc on $\hat \Delta^\pm\setminus (\d \gamma^\pm\times \int D^2)$ so that $\beta^\pm[a^\pm]$ links $G(S^1\times a^\pm)$ zero times; see  Figure \ref{fig:fig-delta-hat-Delta-beta-plus-minus}. Let $\hat \beta^\pm$ be a pushoff of $\beta^\pm$ along  a normal vector field tangent to  $\Delta^\pm$. We see that $\hat \beta^\pm$ and $\hat \beta^\pm$ bound embedded 2-disks $F^\pm$ and $\hat F^\pm$ in $D^3$, respectively, which are disjoint (that is to say, the aforementioned normal vector field extends to a normal vector field of $F^\pm$). Then $\beta^\pm[a^\pm]$ and $\hat \beta^\pm[a^\pm]$ bound the disjoint, embedded 2-disks $F^\pm[a^\pm]$ and $\hat F^\pm[a^\pm]$, respectively, which lie in  $D^3[a^\pm]\setminus G(S^1\times a^\pm)$.  

\begin{figure}[h]
\centering
    \includegraphics[width=.9\linewidth]{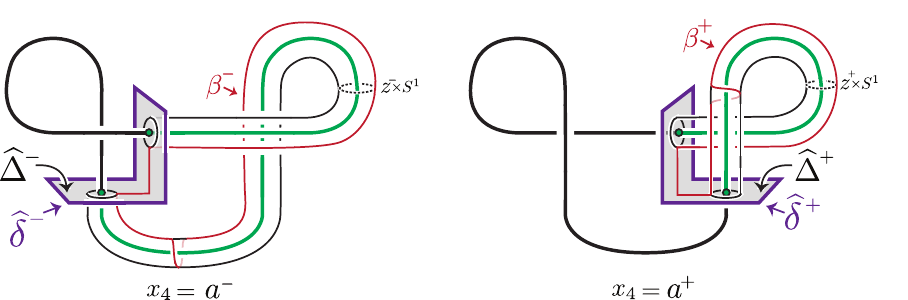}
    \caption{}
    \label{fig:fig-delta-hat-Delta-beta-plus-minus}
\end{figure}

Since $f(S^2_+)$ is disjoint from $\Theta_i(N^+[-1,1])$ and (we may assume) from a collar of the 2-disk $\Theta_i(F^+[a^+])$, observe that $F^+$ may be constructed from a normal pushoff ${\hat E^+}$ of $E^+$ in $D^3$ by attaching a collar so that the intersections between $\Theta_i(F^+[a^+])$ and $f(S^2_+)$ occur entirely on $\Theta_i(\hat E^+[a^+])$. Similarly, since $f(S^2_+)$ is disjoint from $\Theta_i(N^+[-1,1])$ and $\Theta_i(N^-[-1,1])$,  and (we may assume) from a collar of the 2-disk $\Theta_i(F^-[a^-])$, the 2-disk $F^-$ may be  chosen to contain a normal pushoff $\hat V$ of $V$  in $D^3$ and a normal pushoff $\check E^+$ of $E^+$ in $D^3$ so that the intersections between $\Theta_i(F^-[a^+])$ and $f(S^2_+)$ occur entirely on $\Theta_i(\check E^+[a^-])$ and $\Theta_i(\hat V[a^-])$.
\mycomment{ 
    From these observations, by choosing the same whisker in $X_-$ for the disk $\Theta_i(F^+[a^+])$ and its subdisk $\Theta_i(\hat E^+[a^+]\subset \Theta_i(F^+[a^+])$, and the same whisker for the disk $\Theta_i(F^-[a^-]$ and its subdisks $\Theta_i(\check E^+[a^-])$ and $\Theta_i(\hat V[a^-])$,  by Proposition \ref{prop:lambda-restricted} we have
}
From these observations, by Proposition \ref{prop:lambda-restricted} we have (by an appropriate choice of whiskers and orientations in  $X_-$)
\begin{align}\label{eqn:fplus-intersect-disk1}
    \lambda\bigl(f(S^2_+), \Theta_i(F^+[a^+])\bigr) = \lambda\bigl(f(S^2_+), \Theta_i(\hat E^+[a^+])\bigr)
\end{align}
and
\begin{align}\label{eqn:fplus-intersect-disk2}
    \lambda\bigl(f(S^2_+), \Theta_i(F^-[a^-])\bigr) = \lambda\bigl(f(S^2_+), \Theta_i(\check E^+[a^-])\bigr) + \lambda\bigl(f(S^2_+), \Theta_i(\hat V[a^-])\bigr)
\end{align}
in $\Z[s,s^{-1}]$. But, as  $\Theta_i(\hat E^+[a^+])$ and $\Theta_1(\check E^+[a^-])$ are each normal pushoffs in $X_-$ of  $E_i^+\defeq \Theta_i(E^+[0])$, the 2-disk bounded by $\alpha_i^+$, we deduce from Equation \eqref{eqn:fplus-intersect-disk1} that
\begin{align}\label{eqn:fplus-intersect-disk3}
    \lambda\bigl(f(S^2_+), \Theta_i(F^+[a^+])\bigr) = q^{(1)}_i(s)
\end{align}
for some integer Laurent polynomial $q^{(1)}_i(s)\in \Z[s,s^{-1}]$ such that
\begin{align}\label{eq:qi1}
    q_i^{(1)}(1)=f(S^2_+)\cdot E_i^+=\lk(f(S^2_+), \alpha_i^+).
\end{align}
Moreover, if $|f(S^2_+) \cap E_i^+| =1$ for some $i$, then (as we are free to choose the orientation and whisker of $\Theta_i(F^+[a^+])$)  we may take $\widehat q^{(1)}_i(s)=1$. Similarly, as  $\Theta_i(\hat V[a^-])$ is a normal pushoff of $V_i$, from Equation \eqref{eqn:fplus-intersect-disk2} we have (by an appropriate choice of orientations and whiskers in  $X_-$)
\begin{align}\label{eqn:fplus-intersect-disk4}
    \lambda\bigl(f(S^2_+), \Theta_i(F^-[a^-])\bigr) = q^{(1)}_i(s) + \lambda\bigl(f(S^2_+), V_i^c\bigr),
\end{align}
where $V_i^c\subset X_-$ is obtained from $V_i$ (which has boundary on $f(S^2_-)$) by removing a collar in $S^4\setminus f(S^2_+)$.

From the (embedded) punctured torus $\Delta^\pm$, we remove the interior of the annulus bounded by $\beta^\pm\cup \hat \beta^\pm$ and attach the 2-disks $F^\pm\cup \hat F^\pm$. We thus obtain an embedded 2-disk   $\hat A^\pm\subset D^3$ which has boundary  $\hat \delta^\pm$ and is such that $\hat A^\pm[a^\pm]$  lies  in  $D^3[a^\pm]\setminus G(S^1\times a^\pm)$. 
Consequently,  $\Theta_i(\hat A^\pm[a^\pm])$ is an embedded 2-disk in $X_-$ with boundary $\Theta_i(\hat \delta^\pm[a^\pm])$, obtained from the embedded punctured torus $\Theta_i(\Delta^\pm[a^\pm])\subset X_-\setminus f(S^2_+)$ by applying the construction of Lemma \ref{lem:surgery} with the 2-disk $\Theta_i(F^\pm[a^\pm])$ and its normal pushoff $\Theta_i(\hat F^\pm[a^\pm])$.

Now, for an interior point $z^\pm$ on $\gamma^\pm$, the loop $z^\pm\times S^1\subset \gamma^\pm\times S^1$ is dual to $\beta^\pm$ and $\hat \beta^\pm$ on $\Delta^\pm$, and is meridinal to $\hat G$ in $D^3\times D^1$. Thus the loop $\Theta_i(z^\pm\times S^1)$ is dual to $\d\Theta_i(F^\pm[a^\pm])$ and $\d\Theta_i(\hat F^\pm[a^\pm])$ on $\Theta_i(\Delta^\pm[a^\pm])$, and represents $s$ or $s^{-1}$ in $\pi_1(X_-)=\Z\langle s\rangle$. By Lemma \ref{lem:surgery}(i) and Equations \eqref{eqn:fplus-intersect-disk3}-\eqref{eqn:fplus-intersect-disk4}, then, we have (by an appropriate choice of orientation and whisker in  $X_-$)
\begin{align}\label{eq:f-hat-A-i-plus}
    \lambda\bigl(f(S^2_+), \Theta_i(\hat A^\pm[a^\pm])\bigr) = (1-s)q^{(2)}_i(s)
\end{align}
and
\begin{align}\label{eq:f-hat-A-i-minus}
    \lambda\bigl(f(S^2_+), \Theta_i(\hat A^\pm[a^\pm])\bigr) = (1-s)q^{(2)}_i(s) + (1-s)\lambda\bigl(f(S^2_+), V_i^c\bigr)
\end{align}
for some integer Laurent polynomial $q^{(2)}_i(s)\in \Z[s,s^{-1}]$ such that $q_i^{(2)}(1)=q_i^{(1)}(1)$.

We proceed to attach an annulus  to each of the 2-disks $\hat A^\pm[a^\pm]$ and $\hat A^\pm[b^\pm]$ so to obtain 2-disks with which to surger the linking torus $T^\pm$ using the construction of Lemma \ref{lem:surgery}.

In Figure \ref{fig:fig-T-delta-a-b-and-iso-copies}  we have illustrated an oriented circle $\delta_a^\pm$ on $T^\pm\subset N^\pm[-1,1]$ which intersects $D^3[-1]$ and $D^3[1]$ each in an arc, and appears as a pair of points in $D^3[t]$ for $t\in (-1,1)$. We have also illustrated a normal pushoff $\delta_b^\pm$ of $\delta_a^\pm$ on $T^\pm$. Let $\hat T^\pm$ denote the annulus on $T^\pm$ bounded by $\delta_a^\pm\cup \delta_b^\pm$. \mycomment{Let $T_\delta$ denote the annulus bounded by $\delta_{-1}^+\cup \delta_{1}^-$ on $T^+$.} 
Notice that the pair $\delta_a^\pm\cup \delta_b^\pm$ is  isotopic in $N^\pm[-1,1]\setminus \hat G$ to the link $\delta_a^{\pm\prime}\cup\delta_b^{\pm\prime}$  in $N^\pm[-1,1]\setminus \hat G$ which we have also illustrated in Figure \ref{fig:fig-T-delta-a-b-and-iso-copies}. We then see that the annulus $\hat T^\pm$ is homotopic in $N^\pm[-1,1]\setminus \hat G$ to the annulus $\hat \delta^\pm[a^\pm, b^\pm]$ by a homotopy of $S^1\times I$ whose restriction to $S^1\times \{0,1\}$ is a homotopy from $\delta_a^\pm\cup \delta_b^\pm$ to the link $\hat \delta^\pm[a^\pm]\cup \hat \delta^\pm[b^\pm]$  through a sequence of links, except for one singular link where the two components pass through each other. That is, we can find a regular homotopy $K^\pm:(S^1\times I)\times I\to N^\pm[-1,1]\setminus \hat G$ such that
\begin{enumerate}[label=\textup(\arabic*)]
    \item $K_0^\pm(S^1\times I) = \hat T^\pm$, $K_0^\pm(S^1\times 0)=\delta_a^\pm$, $K_0^\pm(S^1\times 1) = \delta_b^\pm$;
    \item $K_1^\pm(S^1\times I) = \hat \delta^\pm[a^\pm,b^\pm]$,  $K_1^\pm(S^1\times 0) = \hat \delta^\pm[a^\pm]$, $K_1^\pm(S^1\times 1) =  \hat \delta^\pm[b^\pm]$; and
    \item $K_t^\pm(S^1\times \{0,1\})$ is a link for all $t\in I$ except at one value $t'$, where $K^\pm_{t'}(S^1\times \{0,1\})$ has precisely with one transverse self-intersection point, arising as an intersection point between  $K^\pm_{t'}(S^1\times 0)$ and $K^\pm_{t'}(S^1\times 1)$. Denote this intersection point by $w^\pm$.
\end{enumerate}
We may further suppose that the image of the homotopy $K^\pm$ lies in a tubular neighborhood ($\approx S^1\times D^3$) of $\delta_a^\pm$ in $N^\pm[-1,1]\setminus \hat G$. Hence, in particular, we may assume that for each $t\in (0,1)$ the image of $K^\pm_t$ is disjoint from each of $\hat A^\pm[a^\pm]$ and $\hat A^\pm[b^\pm]$. 

\begin{figure}%
    \centering
    \subfloat[The torus $T^+$ in $N^+\times D^1$]{%
    \label{fig:fig-T-plus-delta-a-b-and-iso-copies}%
    \includegraphics[width=\linewidth]{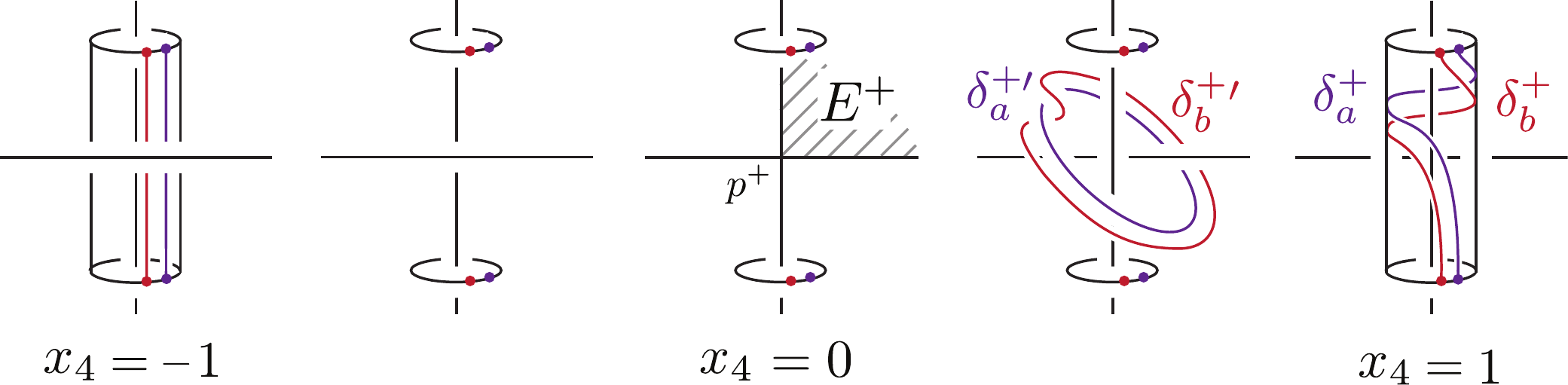}%
    }%
    \vskip 0.5cm
    \subfloat[The torus $T^-$ in $N^-\times D^1$]{%
    \label{fig:fig-T-minus-delta-a-b-and-iso-copies}%
    \includegraphics[width=\linewidth]{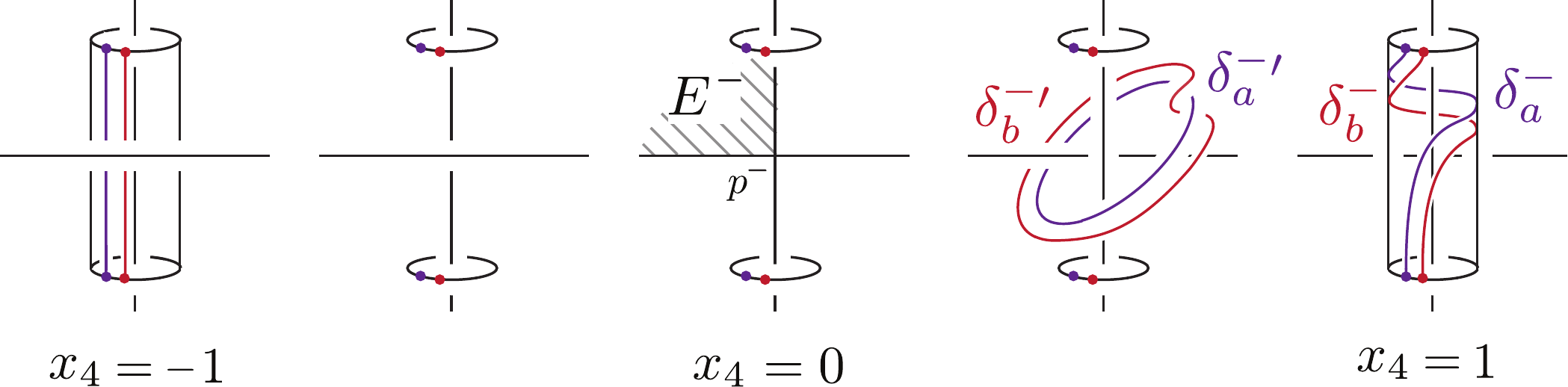}%
    }%
    \caption{}%
    \label{fig:fig-T-delta-a-b-and-iso-copies}%
\end{figure}

Attaching the annuli $K^\pm((S^1\times 0)\times I)$ and $K^\pm((S^1\times 1)\times I)$ to the 2-disks $\hat A^\pm[a^\pm]$ and $\hat A^\pm[b^\pm]$ along $\hat \delta^\pm[a^\pm]=K^\pm(S^1\times 0\times 1)$ and $\hat \delta^\pm[b^\pm]=K^\pm(S^1\times 1\times 1)$, respectively, 
we obtain embedded 2-disks
\begin{align}\label{eq:Omega-a}
    \Omega^\pm_{a} = \hat A^\pm[a^\pm] \mycup{\hat \delta^\pm[a^\pm]}{} K^\pm((S^1\times 0)\times I)
\end{align}
and
\begin{align}\label{eq:Omega-b}
    \Omega^\pm_{b} = \hat A^\pm[b^\pm] \mycup{\hat \delta^\pm[b^\pm]}{} K^\pm((S^1\times 1)\times I)
\end{align}
in $D^4\setminus \hat G$. Observe that their respective boundaries $\d \Omega^\pm_{a} = \delta^\pm_{a}$ and $\d \Omega^\pm_{b} = \delta^\pm_{b}$ lie on $T^\pm$.

By construction,  $\Omega^\pm_a$ and $\Omega^\pm_b$ intersect precisely once, transversely, at the intersection point $w^\pm$ between $K^\pm((S^1\times 0)\times I)$ and $K^\pm((S^1\times 1)\times I)$. Also, since the image of $K^\pm$ lies in $N^\pm[-1,1]$, the intersections between $\Theta_i(\Omega^\pm_a)$ and $f(S^2_+)$ lie precisely on $\Theta_i(\hat A^\pm[a^\pm])$. Consequently, by Proposition \ref{prop:lambda-restricted} we have (after an appropriate choice of orientations and whiskers in  $X_-$)
\begin{align}\label{eq:f-lambda-Omega}
    \lambda\bigl(f(S^2_+),\Theta_i(\Omega^\pm_a)\bigr) = \lambda\bigl(f(S^2_+), \Theta_i(\hat A^\pm[a^\pm])\bigr).
\end{align}

Now, let $A^\pm$ be the 2-sphere in $D^4\setminus \hat G$  obtained by removing from the (embedded) torus $T^\pm$ the interior of the annulus $\hat T^+$ (bounded by $\delta_a^\pm\cup \delta_b^\pm$) and attaching the embedded 2-disks $\Omega^\pm_{a}\cup \Omega^\pm_{b}$. Then $A^\pm$ is immersed and self-transverse in $D^4$,  with a single double point: the point $w^\pm$.

Wishing to apply Lemma \ref{lem:surgery}, define a homotopy from $\Omega_a^\pm$ to $\Omega_b^\pm$ as follows. Let $c:S^1\times I\to D^2$ be a collar with $c(S^1\times 1)=\d D^2$, and let $F_{\hat A^\pm}:D^2\to D^3$ be an embedding with image $\hat A^\pm$.  Since $K^\pm(S^1\times I\times I)$ and $\hat A^\pm[a,b]$ each lie in $D^4\setminus \hat G$,  from Equations  \eqref{eq:Omega-a}-\eqref{eq:Omega-b} it is readily seen that
\begin{align*}
    J^\pm_t(c(x,s)) &= K^\pm(x,t,s)  && \text{ for $(x,s)\in S^1\times I$,}\\
    J^\pm_t(y) &= F_{\hat A^\pm}(y)[a^\pm+t(b^\pm - a^\pm)]  && \text{ for $y\in D^2$,}
\end{align*}
defines a homotopy  $J^\pm:D^2\times I\to D^4\setminus \hat G$ from $\Omega_a^\pm$ to $\Omega_b^\pm$.
Then
\begin{align}\label{eq:A-def}
    A^\pm = (T^\pm \setminus \int \hat T^\pm) \mycup{\delta_a^\pm\cup \delta_b^\pm}{} J^\pm(D^2\times \{0,1\}).
\end{align}
 Now, let $A_i^\pm = \Theta_i(A^\pm)$; then  $A_i^\pm$ is an immersed, self-transverse 2-sphere in $D^3[t_i^-,t_i^+]\setminus f(S^2_-)$ constructed by surgering the embedded torus $T_i^\pm\subset X_-\setminus f(S^2_+)$ using the 2-disks $\Theta_i(\Omega^\pm_{a})\cup \Theta_i(\Omega^\pm_{b})\subset X_-$. Furthermore, observe that a dual curve to $\delta_a^\pm$ and $\delta_b^\pm$ on the torus $T^\pm$ is meridinal to $\hat G$ in $D^4$, so a dual curve to $\Theta_i(\delta_a^\pm)$ and $\Theta_i(\delta_a^\pm)$ on the linking torus $\Theta_i(T^\pm)$ represents $s$ or $s^{-1}$ in $\pi_1(X_-)=\Z\langle s\rangle$. Thus, by Equation \eqref{eq:A-def} and Lemma \ref{lem:surgery}(i) we have  (after an appropriate choice of whiskers and orientations  in $X_-$)
\[
    \lambda\bigl(f(S^2_+), A_i^+\bigr) = (1-s)\lambda\bigl(f(S^2_+), \Theta_i(\Omega_a^\pm)\bigr).
\]
From Equations \eqref{eq:f-hat-A-i-plus} and \eqref{eq:f-lambda-Omega} we therefore have 
\begin{align}\label{eq:f-A-i-plus}
    \lambda\bigl(f(S^2_+), A_i^+\bigr) =(1-s)^2q_i(s) 
\end{align}
and
\begin{align}\label{eq:f-A-i-minus}
    \lambda\bigl(f(S^2_+), A_i^-)\bigr) = (1-s)^2q_i(s) + (1-s)^2\lambda\bigl(f(S^2_+), V_i^c\bigr)
\end{align}
for some $q_i(s)\in \Z[s,s^{-1}]$ such that $q_i(1)=\lk(f(S^2_+), \alpha_i^+)$. As observed earlier, if $|f(S^2_+) \cap E_i^+| =1$ for some $i$, then  we may take $q_i(s)=1$. Moreover, since $\Theta_i(\Omega^\pm_{a})$ and $\Theta_i(\Omega^\pm_{b})$ are embedded and intersect precisely once, we have from Lemma \ref{lem:surgery}(ii) that
\[
    \mu(A_i^\pm) = s \mod 2
\]
in $\Z_2[s]$. Hence $\lambda(A_i^\pm, A_i^{\pm\prime}) = s+s^{-1}\mod 2$.

Finally, by construction, $A_i^\pm$ is homologous to $T_i^\pm$  for each $i$, so by \cite[Lemma 4.3]{me2} \mycomment{Lemma \ref{lem:pi2}}  the immersed 2-spheres $\{A_i^+, A_i^-\}_{i=1}^d$ represent a $\Z[s,s^{-1}]$-basis  for $\pi_2(X_-)$.
\end{proof}

The rest of this section will be devoted to applying Lemma \ref{prop:exchange-all-b-i} to prove the following proposition, which will allow us to surger out the intersections between each $V_i$ and $f(S^2_+)$ (in exchange for intersections with $f(S^2_-))$.

\begin{proposition}\label{prop:FM-constructVis-withu}
Let $f$ be a good link map such that $\sigma_-(f)=0$ and $f(S^2_-)=\hat \U_d$ for some $d\geq 0$.  Then, perhaps after  an ambient isotopy, we may assume that $f(S^2_-)=\hat \U_d$ and the embedded Whitney disks $\{V_i\}_{i=1}^d$ in $S^4$ are framed and satisfy $V_i\cdot f(S^2_+) = 0$ for each $1\leq i\leq d$.
\end{proposition}

The remainder of this section shall be devoted to proving this result. Recall from the beginning of the present section that on the immersed circle $G(S^1\times 0)$ in $D^3[0]$,  the arc $G(D^1_\pm\times 0)$ contains the loop  $\alpha^\pm=G(\hat \alpha^\pm\times0)$ in its interior.

\begin{lemma}\label{lem:permute-dps-second-u}
Let $d\geq 0$. For each $i\in\{1,2,\ldots, d\}$ let $\mu_i\in\{-1,1\}$, and let $\varsigma$ be a permutation of $\{1,2,\ldots, d\}$. There are 4-ball neighborhoods $B_i^\pm$ of $G_i(D^1_\pm\times S^1)$ in $S^4$, $i=1,\ldots, d$, and  an ambient isotopy $\varphi:S^4\times I\to S^4$ such that $\varphi_1(\hat \U_d)=\hat \U_d$, and  for each $1\leq i\leq d$,
\begin{enumerate}[label=\textup{(\roman*)}, ref=\textup{(\roman*)},align=CenterWithParen]
    \item\label{item:permute1} $\varphi_1$ restricts to the identity on $B_i^+$ (so $\varphi_1\circ \alpha_i^+ = \alpha_i^+$),
    \item\label{item:permute2} $\varphi_1$ carries $(B_i^-, G_i(D^1_-\times D^1))$ to $(B_{\varsigma(i)}^-, G_{\varsigma(i)}(D^1_-\times D^1))$ and
    \item\label{item:permute3} {\color{white}{.}}\vspace*{-.6cm}\begin{align*} \varphi_1\circ \alpha_i^- = \begin{cases} \alpha_{\varsigma(i)}^- & \text{ if $\mu_i=1$,} \\[0.3cm] \overline{\alpha_{\varsigma(i)}^-} & \text{ if $\mu_i=-1$.} \end{cases}\end{align*}
\end{enumerate}
\end{lemma}
\begin{proof}
By the construction of $G$, there are disjoint 3-balls $B^3_+$ and $B^3_-$ in $D^3$  such that $B^3_\pm[0]$ is a neighborhood of $ E^\pm\subset D^3[0]$ and such that there is an orientation-preserving diffeomorphism $\Pi^\pm:D^4\to B^3_\pm\times D^1$ carrying the cusp $D^\pm$ of Section \ref{sec:unknotted} to $G(D^1_\pm\times D^1)$, the double point $r^\pm$  to $p^\pm$, and the oriented loop $\theta^\pm$ to $\alpha^\pm$.  \mycomment{ and $\d(B^3_\pm\times D^1) = \d B^3_\pm\times D^1\cup B^3_\pm\times \d I$ precisely along $\d G(D^1_\pm\times D^1) = G(\d D^1_\pm\times D^1)\cup G(D^1_\pm\times \d I)$.}

For each $1\leq i\leq d$, let $\hat \mu_i\in \{-,+\}$ denote the sign of $\mu_i$, let $\Psi_i^\pm$ be the orientation-preserving diffeomorphism given by
\[
    \Psi_i^\pm = \Theta_i\circ \Pi^\pm: D^2\times D^2 \to B^3_\pm\times I_i,
\]
and let $\overline{\Psi_i^\pm}=\Psi_i^\pm \circ \Sigma$, where $\Sigma$ is the orientation-preserving diffeomorphism of $D^2\times D^2$  defined in Section \ref{sec:unknotted} by $\Sigma(x,y)=(-x,y)$. Recalling Equations \eqref{eq:propertiesofSigma-1}-\eqref{eq:propertiesofSigma-2}, observe that $\Psi_i^\pm(D^\pm) = G_i(D^1_\pm\times D^1) = \overline{\Psi_i^\pm}(D^\pm)$,
\begin{align}\label{eq:permute-psi1-alpha}
    \Psi_i^\pm\circ \theta^\pm = \alpha_i^\pm \text{ and } \overline{\Psi_i^\pm} \circ \theta^\pm = \overline{\alpha_i^\pm}.
\end{align}
Furthermore, by construction, $\hat \U_d$ is obtained from the unknotted, embedded 2-sphere $(S^1\times D^1)\cup (D^2\times \{\pm 1\})\subset D^3\times D^1$ by removing its intersections with the 4-balls $B^3_\pm\times I_i = \Psi_i^\pm(D^4)$, and attaching the cusps  $G_i(D^1_\pm\times D^1) = \Psi_i^\pm(D^\pm)$, for  $i=1,2,\ldots, d$. For each $1\leq i\leq d$, let $b_{2i}=\Psi_i^+$, $b_{2i-1} = \Psi_i^-$,  and define a permutation $\rho$ on $\{1,2,\ldots, 2d\}$ by $\rho(2i)=2i$ and $\rho(2i-1)=2\varsigma(i)-1$. Then  Lemma \ref{prop:exchange-all-b-i} yields an ambient isotopy $\varphi:S^4\times I\to S^4$ such that $\varphi_1$ fixes
\begin{align}\label{eqn:permute-ud}
    \hat \U_d\setminus \int \mycup{i=1}{d} (b_{2i}(D^4)\cup b_{2i-1}(D^4)) = \hat \U_d\setminus \int \mycup{i=1}{d} ((B^3_+\cup B^3_-)\times I_i)
\end{align}
set-wise, satisfies
\begin{align}
    \varphi_1\circ \Psi_i^+ &= \varphi_1\circ b_{2i} = b_{2i} = \Psi_i^+ \text{ and} \label{eqn:psi1}\\
    \varphi_1\circ \Psi_i^-(x,y) &= \varphi_1\circ b_{2i-1}(x,y) = b_{2\varsigma(i)-1}(\hat\mu_i\, x, y)=\Psi_{\varsigma(i)}^-(\hat\mu_i\, x, y)\label{eqn:psi2}.
\end{align}
Now,  putting $B_i^\pm=B^\pm\times I_i$, Equation \eqref{eqn:psi1} gives part \ref{item:permute1} of the lemma; Equation \eqref{eqn:psi2} gives part \ref{item:permute2}, and part \ref{item:permute3} follows from Equation \eqref{eq:permute-psi1-alpha} and by noting that Equation \eqref{eqn:psi2} implies $\varphi_1\circ \Psi_i^- =\Psi_{\varsigma(i)}^-$ if $\mu_i=1$ and $\varphi_1\circ \Psi_i^- =\overline{\Psi_{\varsigma(i)}^-}$ if $\mu_i=-1$. Since $\overline{\Psi_i^-}(D^-)= \Psi_i^-(D^-)$, $\varphi_1$ sends $\mycup{i=1}{d} (\Psi_i^+(D^+) \cup \Psi_i^-(D^-))$ to
\[
    \mycup{i=1}{d} (\Psi_i^+(D^+) \cup \Psi_{\varsigma(i)}^-(D^-)) = \mycup{i=1}{d} (\Psi_i^+(D^+) \cup \Psi_{i}^-(D^-)),
\]
so $\varphi_1(\hat \U_d) = \hat \U_d$ by Equation \eqref{eqn:permute-ud}.
\end{proof}

We may now perform an ambient isotopy which carries $f(S^2_-)=\hat \U_d$ back to itself in such a way that the accessory circles $\{\alpha_i^+, \alpha_i^-\}_{i=1}^d$ are rearranged into canceling pairs with respect to their linking numbers with $f(S^2_+)$.

\begin{lemma}\label{lem:linking-same-in-abs-value}
Let $f$ be a \mycomment{good }link map such that $\sigma_-(f)=0$ and  $f(S^2_-)=\hat \U_d$ for some $d\geq 0$. 
Then $f$ is link homotopic (in fact, ambient isotopic) to a link map $g$ such that $g(S^2_-) =  \hat \U_d$ and, for each $1\leq i\leq d$,
\[
    \lk\bigl(\alpha_i^+, g(S^2_+)\bigr) = \lk\bigl(\alpha_i^-, g(S^2_+)\bigr).
\]
\end{lemma}
\begin{proof}
Since $\sigma_-(f)=0$, there is a function $ \mu:\{1,2,\ldots, d\}\to\{-1,1\}$ and a permutation $\varsigma$ on $\{1,2,\ldots, d\}$ such that
\[
    \lk\bigl(\alpha_i^+, f(S^2_+)\bigr) =  \mu_i \cdot \lk\bigl(\alpha_{\varsigma^{-1}(i)}^-, f(S^2_+)\bigr)
\]
for each $1\leq i\leq d$. By Lemma \ref{lem:permute-dps-second-u}, there is an  ambient isotopy $\varphi:S^4\times I\to S^4$ such that  $\varphi_1(\hat \U_d)=\hat \U_d$ and,  for each $1\leq i\leq d$,  $\varphi_1\circ \alpha_i^+ =  \alpha_i^+$, $\varphi_1\circ \alpha_i^- = \alpha_{\varsigma(i)}^-$ if $\mu_i=1$, and
\[
    \varphi_1\circ \alpha_i^- = \overline{\alpha_{\varsigma(i)}^-}
\]
if $\mu_i=-1$. Then, for each $1\leq i\leq d$,
\begin{align*}
    \varphi_1^{-1}(\alpha_i^-) = \begin{cases} \alpha_{\varsigma^{-1}(i)} & \text{ if $\mu_i=1,$}\\ \overline{\alpha_{\varsigma^{-1}(i)}} & \text{ if $\mu_i=-1$}, \end{cases}
\end{align*}
and hence
\begin{align}
    \lk\bigl(\varphi_1^{-1}(\alpha_i^-), f(S^2_+)\bigr) &= \mu_i \cdot \lk\bigl(\alpha_{\varsigma^{-1}(i)}, f(S^2_+)\bigr)\notag\\
     &= \lk\bigl(\alpha_i^+, f(S^2_+)\bigr).
\end{align}
Thus, taking $g=\varphi_1\circ f$, we have
\begin{align*}
    \lk\bigl(\alpha_i^-, g(S^2_+)\bigr) &=  \lk\bigl(\varphi_1(\varphi_1^{-1}(\alpha_i^-)), \varphi_1(f(S^2_+))\bigr)\\
    &= \lk\bigl(\varphi_1^{-1}(\alpha_i^-), f(S^2_+)\bigr)\\
    &=  \lk\bigl(\alpha_i^+, f(S^2_+)\bigr)\\
    &= \lk\bigl(\varphi_1(\alpha_i^+), \varphi_1(f(S^2_+))\bigr)\\
    &= \lk\bigl(\alpha_i^+, g(S^2_+)\bigr).\qedhere
\end{align*}
\end{proof}

Having established a means to permute the accessory circles of $\hat \U_d$ in a prescribed way, we may now complete the proof of Proposition \ref{prop:FM-constructVis-withu}.

\begin{proof}[Proof of Proposition \ref{prop:FM-constructVis-withu}]
By Lemma \ref{lem:linking-same-in-abs-value} we may assume, after an ambient isotopy, that  $f(S^2_-)=\hat \U_d$ and the accessory circles $\{\alpha_i^+, \alpha_i^-\}_{i=1}^d$ on $f(S^2_-)$ satisfy $\lk\bigl(\alpha_i^+, f(S^2_+)\bigr) = \lk\bigl(\alpha_i^-, f(S^2_+)\bigr)$ for each $1\leq i\leq d$.
Recall the notation of Figure \ref{fig:circle} and that we  let $\hat \eta_x$ and $\hat \eta_y$ denote disjoint neighborhoods of $\eta_x$ and $\eta_y$, respectively, on the circle $S^1$. For each $1\leq i\leq d$, $G_i(\hat \eta_x\times [-\tfrac{1}{2}, \tfrac{1}{2}])$ and $G_i(\hat \eta_y\times [-\tfrac{1}{2}, \tfrac{1}{2}])$ are embedded 2-disk neighborhoods of $\{p_i^+,p_i^-\}$ on $f(S^2_-)$ which  intersect  precisely at these two points, and the accessory circles $\{\alpha_i^+, \alpha_i^-\}$  leave along $G_i(\hat \eta_x\times [-\tfrac{1}{2}, \tfrac{1}{2}])$ and return along $G_i(\hat \eta_y\times [-\tfrac{1}{2}, \tfrac{1}{2}])$. Thus, as the arc $G_i(\eta_x\times 0)\subset G_i(\hat \eta_x\times [-\tfrac{1}{2}, \tfrac{1}{2}])$ runs from $p^+$ to $p^-$, and the arc $G_i(\eta_y\times 0)\subset G_i(\hat \eta_y\times [-\tfrac{1}{2}, \tfrac{1}{2}])$ runs from $p^+$ to $p^-$, by Lemma \ref{lem:linking-number-whitney-circle} and Equation \eqref{eq:bdy-V-i} we have
\begin{align*}
    \bigl|\lk(\d V_i, f(S^2_+))\bigr| &= \bigl|\lk(G_i(\eta_x\times 0)\cup G_i(\eta_y\times 0), f(S^2_+))\bigr|\\
     &= \bigl|\lk(\alpha_i^+, f(S^2_+)) - \lk(\alpha_i^-, f(S^2_+))\bigr| \\
     &= 0.\qedhere
\end{align*}
\end{proof}

\subsection{Whitney disks in $X_+$}\label{sec:omega}



Referring to the notation of Proposition \ref{prop:FM-constructVis-withu} and Proposition \ref{lem:2-spheres}, we next show that by altering the interiors of the 2-disks $\{V_i\}_i$ so to exchange their intersections with $f(S^2_+)$ for intersections with  $f(S^2_-)$, we are able to compute $\omega_-$ as follows.

\begin{proposition}\label{prop:relate-omega-Vi}
For each $1\leq i\leq d$, the pair $\{p_i^+, p_i^-\}$ of double points of $f(S^2_-)\subset X_+$ may be equipped with a framed, immersed Whitney disk $W_i$ in $X_+$ such that $\d W_i=\d V_i$. Furthermore, there are integer Laurent polynomials $\{u_i(s)\}_{i=1}^d$ such that
\[
    \omega_-(f) = \mysum{i: \text{ $n_i$ even}}{} u_i(1) \mod 2
\]
and for each $1\leq i\leq d$,
\[
    \lambda(f(S^2_+), V_i^c) = (1+s)u_i(s)  \mod 2.
\]
\end{proposition}

Define a ring homomorphism $\varphi:\Z[s,s^{-1}] \to \Z_2$ by
\[
    \Z[s,s^{-1}]\xrightarrow{\d} \Z[s,s^{-1}]\xrightarrow{s\,\mapsto 1}\Z \xrightarrow{\text{mod } 2} \Z_2,
\]
where $\d$ is the formal derivative defined by setting $\d(s^n)=ns^{n-1}$ (for $n\in \Z$) and extending by linearity. Recall from Section \ref{sec:prelims} that we use $\equiv$ to denote equivalence modulo $2$, and for an integer Laurent polynomial $g(s)\in \Z[s,s^{-1}]$ we write $\overline {g(s)} = g(s^{-1})$. The following properties of $\varphi$ are readily verified.
\begin{lemma}\label{lem:varphi-props}
If $g(s)\in \Z[s,s^{-1}]$, then
\begin{enumerate}[label=\textup{(\roman*)}, ref=\textup{(\roman*)},align=CenterWithParen]
    \item $\varphi(\overline{g(s)}) \equiv \varphi(g(s))$,
    \item $\varphi(s\cdot g(s)) \equiv g(1) + \varphi(g(s))$, and
    \item $\varphi((1+s^n)g(s)) \equiv n\cdot g(1)$.
\end{enumerate}
\end{lemma}

Let $1\leq i\leq d$. Since $f(S^2_+)\cdot V_i=0$ the intersections between $f(S^2_+)$ and $\int V_i^c$ (which may be assumed transverse after a small homotopy of $f_+$) may be decomposed into pairs of opposite sign $\{x_i^j, y_i^j\}_{j=1}^{J_i}$ for some $J_i\geq 0$ (for any choice of orientation of $S^2_+$ and $V_i^c$). For each $1\leq j\leq J_i$, choose a simple path $\alpha_i^j$ on $f(S^2_+)$ from $x_i^j$ to $y_i^j$ whose interior is disjoint from $\mycup{k=1}{d} V_k$,  let $\beta_i^j$ be a simple path in $\int V_i^c$ from $y_i^j$ to $x_i^j$ whose interior misses $f(S^2_+)$, and let $\rho_i^j=\alpha_i^j\cup \beta_i^j$. The resulting collection of loops $\{\rho_i^j\}_{j=1}^{J_i}$ in $X_-$ may be chosen to be mutually disjoint. For each $1\leq j\leq J_i$ define the $\Z_2$-integer
\[
    m_i^j = \lk(f(S^2_-), \rho_i^j) \mod 2.
\]
Note that $m_i^j$ is well-defined because $f(S^2_-)$ and $V_i^c$ are simply-connected (c.f. Proposition \ref{prop:lambda-product}).
\begin{lemma}\label{lem:m-i-j-new}
There are integer Laurent polynomials $\{u_i(s)\}_{i=1}^d$  in $\Z[s,s^{-1}]$ such that for each $1\leq i\leq d$ we have
\[
    \lambda(f(S^2_+), V_i^c) \equiv (1+s)u_i(s)
\]
and $u_i(1) \equiv \sum_{j=1}^{J_i} m_i^j$.
\end{lemma}
\begin{proof}
Choose whiskers connecting $f(S^2_+)$ and $\{V_i^c\}_{i=1}^d$ to the basepoint of $X_-$. 
Let $1\leq i\leq d$ and $1\leq j\leq J_i$. Since $\rho_i^j$ is a loop in $X_-$ that runs from $x_i^j$ to $y_i^j$ along $f(S^2_+)$, and back to $x_i^j$ along $V_i^c$, by Proposition \ref{prop:lambda-product} we have
\[
    \lambda(f(S^2_+), V_i^c)[x_i^j]\cdot (\lambda(f(S^2_+), V_i^c)[y_i^j])^{-1} = s^{\hat m_i^j} \in \pi_1(X_-),
\]
where $\hat m_i^j$ is an integer  such that $\hat m_i^j \equiv \lk(f(S^2_-), \rho_i^j) \equiv  m_i^j.$  Thus
\[
   \lambda(f(S^2_+), V_i^c)[x_i^j] = s^{\hat m_i^j}\cdot \lambda(f(S^2_+), V_i^c)[y_i^j],
\]
so the mod $2$ contribution to $\lambda(f(S^2_+), V_i^c)$ due to the pair of intersections $\{x_i^j, y_i^j\}$ is
\[
    \lambda(f(S^2_+), V_i^c)[x_i^j] + \lambda(f(S^2_+), V_i^c)[y_i^j] \equiv  (1+s^{\hat m_i^j})s^{l_i^j},
\]
for some $l_i^j\in \Z$.
Choose $u_i^j(s)\in \Z[s,s^{-1}]$ such that
\begin{align*}
    (1+s)u_i^j(s) \equiv (1+s^{\hat m_i^j})s^{l_i^j};
\end{align*}
applying $\varphi$ to both sides (c.f. Lemma \ref{lem:varphi-props}) yields
\begin{align}\label{eq:u-i-j-m-i-j-new}
    u_i^j(1)\equiv m_i^j.
\end{align}
Summing over all such pairs $\{x_i^j, y_i^j\}_{j=1}^{J_i}$ we have
\[
    \lambda(f(S^2_+), V_i^c) \equiv \mysum{j=1}{J_i} (1+s)u_i^j(s)\equiv  (1+s)u_i(s),
\]
where $u_i(s) = \sum_{j=1}^{J_i} u_i^j(s)$ satisfies $u_i(1) \equiv \sum_{j=1}^{J_i} m_i^j$ by Equation \eqref{eq:u-i-j-m-i-j-new}.
\end{proof}

Let $1\leq i\leq d$. We now  remove the intersections between $V_i$ and $f(S^2_+)$ by surgering $V_i$ along the paths $\{\alpha_i^j\}_{j=1}^{J_i}$, obtaining an embedded $J_i$-genus, once-punctured surface $\hat V_i$ in $X_+$ which has interior in $X_-$ and coincides with $V_i$ near the boundary.

\begin{figure}[H]
\centering
    \includegraphics[width=.45\textwidth]{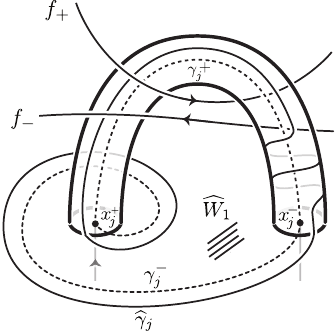}
    \caption{}
    \label{fig:WD-handle}
\end{figure}

Since $f(S^2_+)$ is transverse to $V_i$, for each $1\leq j\leq J_i$ the restriction of a tubular neighborhood of $f(S^2_+)$ to the arc $\alpha_i^j$   may be identified with a 3-ball $h_i^j:D^1\times D^2\to X_-$ such that $h_i^j(D^1\times 0)=\alpha_i^j$ and $h_i^j(D^1\times D^2)$ intersects $V_i$ in two embedded 2-disks $h_i^j(1\times D^2)$ and $h_i^j(-1\times D^2)$ neighborhoods of  $x_i^j$ and $y_i^j$ in $V_i$, respectively. Attaching handles to $V_i$ along the arcs $\alpha_i^j$, $j=1,\ldots, J_i$, yields the surface
\[
    \hat V_i = \bigl[V_i \setminus  \mycup{j=1}{J_i} \int h_i^j(\d D^1\times D^2) \bigr] \mycup{j=1}{J_i} \mycup{h_i^j(\d D^1\times \d D^2)}{} h_i^j(D^1\times \d D^2),
\]%
which is disjoint from both $f(S^2_+)$ and $f(S^2_-)$. See Figure \ref{fig:WD-handle}.

Now, for each $1\leq j\leq J_i$ we may assume that $\beta_i^j$ intersects $h_i^j(1\times \d D^2)$ and $h_i^j(-1\times \d D^2)$ exactly once, at points $\hat x_i^j$ and $\hat y_i^j$, respectively. Let $\myp{\beta}_i^j$ be the subarc of $\beta_i^j$ on $\hat V_i$ running from $\hat y_i^j$ to $\hat x_i^j$, let ${\myp{\alpha}_i}^j$ be a path on $h_j(D^1\times \d D^2)$ connecting $\hat x_i^j$ to $\hat y_i^j$, and put $\myp{\rho}_i^j = \myp{\alpha}_i^j\cup \myp{\beta}_i^j$. By band-summing $\myp{\alpha}_i^j$ with meridinal circles  of $f(S^2_+)$  of the form $h_j(z\times \d D^2)$ (for a point $z$ in $\int D^1$) if necessary, we may assume that $\lk(\myp{\rho}_i^j, f(S^2_+)) = 0$ (see Figure \ref{fig:WD-handle}). Hence, as $\pi_1(X_+)$ is abelian, there is an immersed 2-disk $\myp{Q}_i^j$ in $X_+$ bound by $\myp{\rho}_i^j$. We may further assume that $\myp{Q}_i^j$ misses a collar of $\d V_i =\d\hat V_i$ and is transverse to $f(S^2_-)$.

\begin{figure}[H]
\centering
    \includegraphics[width=.4\textwidth]{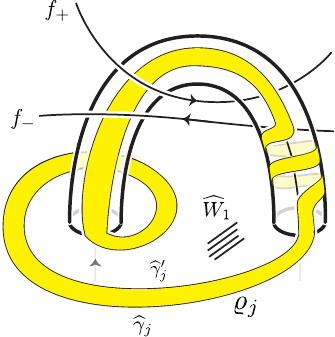}
    \caption{}
    \label{fig:WD-handle-surg}
\end{figure}

By boundary twisting $\myp{Q}_i^j$ along $\myp{\beta}_i^j$ (and so introducing intersections between the interior of $\myp{Q}_i^j$ and $\hat V_i$) if necessary we may further assume that a normal section of $\myp{\rho}_i^j$ that is tangential to $\hat V_i$ extends to a normal section of $\myp{Q}_i^j$ in $X_+$. Hence there is a normal pushoff $\mym{Q}_i^j$ of $\myp{Q}_i^j$ and an annulus $\varrho_i^j$ on $\hat V_i$ with boundary $\d\varrho_i^j=\d \myp{Q}_i^j\cup \d\mym{Q}_i^j$ (see Figure \ref{fig:WD-handle-surg}). Iterating the construction of Lemma \ref{lem:surgery} we may then surger $\hat V_i$ along $\myp{\rho}_i^j$, using $\myp{Q}_i^j$ and its pushoff $\mym{Q}_i^j$, for all $1\leq j\leq J_i$, to obtain an \emph{immersed} 2-disk $W_i$ in $X_+$ such that the framing of $V_i$ (which agrees with $W_i$ near the boundary) along its boundary extends over $W_i$. But $V_i$ is a framed Whitney disk for $f(S^2_-)$ in $S^4$, so $W_i\subset X_+$ is a framed Whitney disk for $f(S^2_-)$ in $X_+\subset S^4$. That is,
\[
    W_i = (\hat V_i \setminus \int \mycup{j=1}{J_i} \,\varrho_i^j) \mycup{j=1}{J_i} \mycup{\d\varrho_i^j}{} (\myp{Q}_i^j\cup \mym{Q}_i^j)
\]
is a framed, immersed Whitney disk for the immersion $f_-:S^2\to X_+$. Let $W_i^c$ denote the complement in $W_i$ of a half-open collar it shares with $V_i$ so that $\d W_i^c = \d V_i^c$ and $f(S^2_-)$ intersects $W_i^c$ in its interior.

The first step in relating $\omega_-(f)$ to the intersections between $f(S^2_+)$ and the $V_i$'s is the following lemma.

\begin{lemma}\label{lem:w-i-n-i-even-new}
The contribution to $\omega_-(f)$ due to intersections between $f(S^2_-)$ and the interior of $W_i$ is
\begin{align*}
    {\cL}^-(W_i) \equiv \begin{cases} 0 & \text{ if $n_i$ is odd,}\\ \mysum{j=1}{J_i} m_i^j & \text{ if $n_i$ is even.}
    \end{cases}
\end{align*}
\end{lemma}

\begin{proof}
Referring to the constructions preceding the lemma statement, since $\int \hat V_i\subset X_-$, the only intersections between $\int W_i$ and $f(S^2_-)$ lie on the immersed 2-disks $\{\myp{Q}_i^j, \mym{Q}_i^j\}_{j=1}^{J_i}$. Indeed, since $\myp{Q}_i^j$ is the pushoff of $\mym{Q}_i^j$ along a section of its normal bundle that is tangent to the annulus $\varrho_i^j$ on  $\hat V_i$, there is an immersion of a 3-ball
\[
    H_i^j: D^2\times I \to X_+
\]
such that $H_i^j(D^2\times 0) = \myp{Q}_i^j$, $H_i^j(D^2\times 1) = \mym{Q}_i^j$ and $H_i^j(S^1\times I)=\varrho_i^j$. Furthermore, since $\myp{Q}_i^j$ is transverse to $f(S^2_+)$ we may assume that if we let
\[
    K_i^j = \#\{f(S^2_-)\cap \myp{Q}_i^j\},
\]
then there are distinct points $x_k\in\int D^2$, $1\leq k\leq K_i^j$, such that $f(S^2_-)$ intersects $H_i^j(D^2\times I)$ precisely along the arcs $\{H_i^j(x_k\times I)\}_k$. Whence the intersections between $f(S^2_-)$ and $\int W_i$ consist precisely of pairs
\begin{align*}
    \myp{x}_i^{j,k} &:= H_i^j(x_k\times 0)\subset \myp{Q}_i^j \text{ and}\\
    \mym{x}_i^{j,k} &:= H_i^j(x_k\times 1)\subset \mym{Q}_i^j,
\end{align*}
for $1\leq k\leq K_i^j$. Thus, in particular, if $n_i$ is odd then from Remark \ref{rem:omega-compute} we have $\cL^-(W_i)=0$.

Suppose now that $n_i$ is even. Note that since the loop $\myp{\rho}_i^j$ on $\hat V_i$ is freely homotopic in $X_-$ to $\rho_i^j$ (to see this, collapse $h_i^j$ onto its core $\alpha_i^j\subset f(S^2_+)$ in $X_-$) we have $\lk(f(S^2_+), \myp{\rho}_i^j)) \equiv m_i^j$, 
so
\begin{align}\label{eq:f-plus-W-i-K-i-j-new}
    K_i^j  \equiv f(S^2_-)\cdot \myp{Q}_i^j  \equiv m_i^j.
\end{align}
We may arrange that there are points $z\in \int D^1$ and $d\in S^1$ so that the meridinal circle $h_i^j(z\times S^1)$ of $f(S^2_+)$ on $\hat V_i$ intersects $\varrho_i^j$ along the arc $H_i^j(d\times I)$.  Let $\iota_I$ denote the interval $[0,1]$ oriented from $0$ to $1$, and let ${\zeta}_i^j$ be the path on $h_i^j(z\times S^1) \setminus \int \varrho_i^j$ that runs from $H_i^j(d\times 0)$\mycomment{$\myp{d}_i^j$} to $H_i^j(d\times 1)$\mycomment{$\mym{d}_i^j$}. Then the loop based at $H_i^j(d\times 1)$\mycomment{$\mym{d}_i^j$} and spanning $h_i^j(z\times S^1)$ given by
\[
    \eta_i^j:= H_i^j(d\times \overline{\iota_I})\ast \zeta_i^j
\]
is a meridinal loop of $f(S^2_+)$, so $\lk(f(S^2_+), \eta_i^j) \equiv 1$.

\renewcommand{\qedsymbol}{$\blacksquare$}
\begin{claim*}\label{claim:omega-contribution}
For fixed $i,j$: if $n_i$ is even then for each $1\leq k\leq K_i^j$ the contribution to ${\cL}^-(W_i)$ due to the pair $\{\myp{x}_i^{j,k}, \mym{x}_i^{j,k}\}$ is
\[
    {\cL}^-_i(\myp{x}_i^{j,k}) + {\cL}^-_i(\mym{x}_i^{j,k}) \equiv 1.
\]
\end{claim*}
\begin{proof}[Proof of claim]
Let $\gamma_i^{j,k}$ be a path in $D^2$ connecting $x_k$ to $d$. Then
\begin{align*}
     \beta_i^{j,k}:=H_i^j(x_k\times \overline{\iota_I})\ast H_i^j(\gamma_i^{j,k}\times 0) \ast \zeta_i^j\ast H_i^j(\overline{\gamma_i^{j,k}}\times 1)
\end{align*}
is a loop that runs from $\mym{x}_i^{j,k}$ to $\myp{x}_i^{j,k}$ along $H_i^j(x_k\times I)\subset f(S^2_-)$ and then back to $\mym{x}_i^{j,k}$ along $W_i^c$, so by Remarks \ref{rem:omega-compute} and \ref{rem:omega-product} we have
\begin{align*}\label{eq:beta-i-j-m-i-j}
    {\cL}^-_i(\myp{x}_i^{j,k}) + {\cL}^-_i(\mym{x}_i^{j,k}) &\equiv m_i(\myp{x}_i^{j,k}) + m_i(\mym{x}_i^{j,k})\\ &\equiv \lk(f(S^2_+), \beta_i^{j,k}).
\end{align*} 
Now, the loop $\beta_i^{j,k}$ is  homotopic in $X_+$ to the loop
\[
    H_i^j(x_k\times \overline{\iota_I})\ast H_i^j(\gamma_i^{j,k}\times 0) \ast H_i^j(d\times \iota_I)\ast \eta_i^j \ast H_i^j(\overline{\gamma_i^{j,k}}\times 1),
\]
but the loop
\[
    H_i^j(x_k\times \overline{\iota_I})\ast H_i^j(\gamma_i^{j,k}\times 0) \ast H_i^j(d\times \iota_I) \ast H_i^j(\overline{\gamma_i^{j,k}}\times 1)
\]
bounds the 2-disk $H_i^j(\gamma_i^{j,k}\times I)\subset X_+$. Thus $\beta_i^{j,k}$ is homotopic in $X_+$ to
\[
    H_i^j(\gamma_i^{j,k}\times 1)\ast \eta_i^j\ast  H_i^j(\overline{\gamma_i^{j,k}}\times 1),
\]
and so $\lk(f(S^2_+), \beta_i^{j,k}) \equiv \lk(f(S^2_+), \eta_i^j) \equiv 1.$
\end{proof}
\renewcommand{\qedsymbol}{$\square$}

Applying the claim to all such pairs of intersections $\{\myp{x}_i^{j,k}, \mym{x}_i^{j,k}\}_{k=1}^{K_i^j}$ between $f(S^2_+)$ and  $\myp{Q}_i^j \cup \mym{Q}_i^j$, over all $1\leq j\leq J_i$, yields the total contribution
\begin{align*}
    \cL^-(W_i) &\equiv  \mysum{x\in f(S^2_-)\cap \,\int W_i}{} {\cL}^-_i(x) \\
    &\equiv \mysum{\;j=1\;}{J_i} \mysum{x\,\in \,\myp{Q}_i^j \cup\, \mym{Q}_i^j}{} {\cL}^-_i(x)\\
    &\equiv \mysum{\;j=1\;}{J_i} \mysum{k=1}{K_i^j} \bigl({\cL}^-_i(\myp{x}_i^{j,k}) + {\cL}^-_i(\mym{x}_i^{j,k})\bigr)\\
    &\equiv \mysum{\;j=1\;}{J_i} K_i^j\\
    &\equiv \mysum{\;j=1\;}{J_i} m_i^j,
\end{align*}
where the last equality is by Equation \eqref{eq:f-plus-W-i-K-i-j-new}. This completes the proof of Lemma \ref{lem:w-i-n-i-even-new}.
\end{proof}

\noindent Applying the lemma, we have
\begin{align*}
    \omega_-(f) &\equiv \mysum{i: \, n_i \text{ odd}}{} \cL^-(W_i) + \mysum{i: \, n_i \text{ even}}{} \cL^-(W_i)\\
    &\equiv \mysum{i: \, n_i \text{ even}}{}\,  \mysum{\;j=1\;}{J_i} m_i^j,
\end{align*}
where $1\leq i\leq d$. Proposition \ref{prop:relate-omega-Vi} now follows from Lemma \ref{lem:m-i-j-new}. \hfill $\square$

\subsection{Relating $\sigma_+$ and $\omega_-$}\label{sec:put-together}\label{sec:main3}

We now bring Kirk's invariant $\sigma_+(f)$ into the picture by noting its relationship with the homotopy class of $f(S^2_+)$ as an element of $\pi_2(X_-)$.

Referring to Proposition \ref{lem:2-spheres}, since $\pi_2(X_-)$ is generated as a $\Z\pi_1(X_-)=\Z[s,s^{-1}]$-module by the 2-spheres $\{A_i^+, A_i^-\}_{i=1}^d$, there are integer Laurent polynomials $\{c_i^+(s), c_i^-(s)\}_{i=1}^d$ such that, as a (whiskered) element of $\pi_2(X_-)$, $f(S^2_+)$ is given by
\begin{align}\label{eq:final-step-1}
f(S^2_+) = \sum_{i=1}^d \,c_i^+(s) A_i^+ + c_i^-(s) A_i^-.
\end{align}
By the sesquilinearity of the intersection form $\lambda(\cdot, \cdot)$ we have from Proposition \ref{lem:2-spheres} that
\begin{align}
    \lambda(f(S^2_+), f(S^2_+)) &\equiv \mysum{i=1}{d} \, c_i^+(s)\overline{c_i^+(s)} \lambda(A_i^+, A_i^+) + c_i^-(s)\overline{c_i^-(s)} \lambda(A_i^-, A_i^-)\notag \\
    &\equiv (s+s^{-1})\mysum{i=1}{d} \bigl[c_i^+(s)\overline{c_i^+(s)} + c_i^-(s)\overline{c_i^-(s)}\bigr].\label{eq:final-step-2}
\end{align}
In \cite{Ki1}, Kirk showed that $\sigma$ has the following image.
\begin{proposition}[\hskip-0.03cm\cite{Ki1}]\label{eq:PK-image}
If $g$ is a link map, then
\[
    \sigma_+(g) + \sigma_-(g) = a_0+ \sum_{n=2}^m a_n(n^2s -s^n)
\]
for some integer $m\geq 0$ and integers $a_0, a_2,a_3,\ldots, a_m$.
\end{proposition}

Now, since $\sigma_-(f)=0$ and  $f$ is a good link map, by Proposition \ref{prop:sigma-mu} we have
\begin{align}
    \lambda(f(S^2_+),f(S^2_+)) &\equiv \sigma_+(f) + \overline{\sigma_+(f)}\notag \\
    &\equiv \sum_{n=2}^m a_n\bigl[s^n + s^{-n} + n(s + s^{-1})\bigr] \label{eq:lambda-n}
\end{align}
for some integers $a_2,\ldots, a_m$. The following observation about the terms in the right-hand side of this equation will be useful in performing some arithmetic in $\Z_2[s,s^{-1}]$.
\begin{lemma}\label{lem:lambda-n}
Let $n\geq 2$ be an integer. Then
\[
    s^n + s^{-n} + n(s + s^{-1}) \equiv (1+s)^4 r_n(s)
\]
for some integer Laurent polynomial $r_n(s)\in \Z[s,s^{-1}]$  such that
\begin{align*}
    r_n(1) = \begin{cases} \frac{n}{2} & \text{ if $n$ is even,}\\
    0 & \text{ if $n$ is odd.}
    \end{cases}
\end{align*}
\end{lemma}

\begin{proof}
If $n=2k$ for some $k\geq 1$, then modulo $2$ we have
\begin{align*}
    s^n + s^{-n} &\equiv s^{-2k}(1 + s^{4k})\\
    &\equiv s^{-2k}(1+s^4)(\underbrace{(s^4)^{k-1}+(s^4)^{k-2}+\ldots + (s^4)^{1} + 1}_{\text{$k$ terms}}).
\end{align*}
On the other hand, if $n=2k+1$ for some $k\geq 1$, by direct expansion one readily verifies that, modulo $2$,
\begin{align*}
    s^n+s^{-1}+s^{-n}+s &\equiv (s+s^{-1})(s^{2k} + s^{2k-2}+\ldots + s^2+1) \\
    & \hspace*{.5cm} + (s+s^{-1})(s^{-2k} + s^{-2k+2}+\ldots + s^{-2}+1)\\
   &\equiv (s+s^{-1})\mysum{l\equiv1}{k} (s^{2l} + s^{-2l})\\
    &\equiv  s^{-1}(1+s^2)\mysum{m\equiv1}{k} s^{-2m} (s^{4m} + 1)\\
    &\equiv s^{-1}(1+s)^2\mysum{m\equiv1}{k} s^{-2m}(s^{m} + 1)^4\\
    &\equiv (1+s)^6 \, \hat r_n(s)
\end{align*}
for some $\hat r_n(s)\in \Z[s,s^{-1}]$.
\end{proof}

\mycomment{
So, as $\Z_2[s,s^{-1}]$ is a domain, \eqref{eq:final-step-2} implies that
\begin{align}\label{eq:c-i}
    \mysum{i=1}{d} c_i^\pm(s)\overline{c_i^\pm(s)} \equiv 0.
\end{align}
}
Now, from Equation \eqref{eq:final-step-1}, Proposition \ref{lem:2-spheres}(ii) and the sesquilinearity of $\lambda(\cdot,\cdot)$, we have
\[
    \lambda(f(S^2_+), A_i^\pm) \equiv c_i^\pm(s) \lambda_2(A_i^\pm, A_i^\pm) \equiv c_i^\pm(s)(s+s^{-1}).
\]
Comparing  with Proposition \ref{lem:2-spheres}(iii),(iv) and Proposition \ref{prop:relate-omega-Vi}, we see that there are integer Laurent polynomials $\{q_i(s), u_i(s)\}_{i=1}^d$ such that
\begin{align}\label{eq:final-step-omega}
    \omega_-(f) \equiv \mysum{i:\text{ $n_i$ even}}{} u_i(1)
\end{align}
and for each $1\leq i\leq d$,
\begin{align*}
    c_i^+(s) &\equiv q_i(s), \\
    c_i^-(s) &\equiv q_i(s) + (1+s)u_i(s),
\end{align*}
where $q_i(1)\equiv n_i$. Thus Equation \eqref{eq:final-step-2} becomes
\begin{align*}
   \;\;\lambda&(f(S^2_+), f(S^2_+))  \\
   \;\;\;\;\;\;&\equiv (s+s^{-1})\mysum{i=1}{d}\bigl[q_i(s)\overline{q_i(s)}
   + (q_i(s) + (1+s)u_i(s))\overline{(q_i(s) + (1+s)u_i(s))}\bigr] \\
     \;\;\;\;\;\;&\equiv (s+s^{-1})\mysum{i=1}{d}\bigl[(1+s)u_i(s)\overline{q_i(s)} + (1+s^{-1})q_i(s)\overline{u_i(s)}  \\
     & \quad\quad\quad +  (1+s)(1+s^{-1})u_i(s)\overline{u_i(s)}\bigr].
\end{align*}
Comparing with Equation \eqref{eq:lambda-n} and applying Lemma \ref{lem:lambda-n} we then have
\begin{align}\label{eq:q-i-u-i-r-n}
    (s+s^{-1})\mysum{i=1}{d}&\bigl[(1+s)u_i(s)\overline{q_i(s)} + (1+s^{-1})q_i(s)\overline{u_i(s)} \\
    & \quad\quad +  (1+s)(1+s^{-1})u_i(s)\overline{u_i(s)}\bigr]\notag\\
    &\equiv (1+s)^4 \sum_{n=2}^k a_n r_n(s) \notag\\
    &\equiv (s+s^{-1})(1+s^{-1})^2 \sum_{n=2}^k a_n \hat r_n(s)
\end{align}
for some integer Laurent polynomials $\{r_n(s), \hat r_n(s)\}_{n=2}^k$ (here, $\hat r_n(s)=s^3r_n(s)$) such that $\hat r_n(1) = r_n(1) = n/2$ if $n$ is even, and $\hat r_n(1)=0$ if $n$ is odd. Since $\Z_2[s,s^{-1}]$ is an integral domain, we may divide both sides of Equation \eqref{eq:q-i-u-i-r-n} by $(s+s^{-1})(1+s^{-1})$ to obtain
\begin{align}\label{eq:summation-2}
      (1+s^{-1}) \sum_{n=2}^k a_n \hat r_n(s) \equiv \mysum{i=1}{d}\, u_i(s)\overline{q_i(s)}s + q_i(s)\overline{u_i(s)} +  (1+s)u_i(s)\overline{u_i(s)}.
\end{align}
Applying the homomorphism $\varphi$ of Lemma \ref{lem:varphi-props} to both sides of Equation \eqref{eq:summation-2} then yields  the following equality in $\Z_2$:
\begin{align*}
     \sum_{n=2}^k a_n \hat r_n(1) &\equiv \mysum{i=1}{d}\bigl[\overline{q_i(1)}u_i(1)+\varphi(u_i(s)\overline{q_i(s)}) + \varphi(q_i(s)\overline{u_i(s)}) + u_i(1)\overline{u_i(1)}\bigr]\\
    &\equiv \mysum{i=1}{d}\, q_i(1)u_i(1) + u_i(1) \\
    &\equiv \mysum{i=1}{d}\, u_i(1)(n_i + 1)\\
    &\equiv \mysum{i: \text{ $n_i$ is even}}{}\, u_i(1).
\end{align*}
Thus, as $\hat r_n(1)\equiv 1$ if and only if $n$ is even and $n/2\equiv  1$ (i.e., $n=2$ mod $4$), from Equation \eqref{eq:final-step-omega} we have
\[
    \omega_-(f) \equiv \sum_n \{a_n:  n = 2 \mod 4\},
\]
completing the proof of Theorem \ref{thm:mainresult}.\qed 


\setcounter{section}{0}
\renewcommand*{\thesection}{\Alph{section}}
\section{Appendix: Proof of Lemma \ref{prop:exchange-all-b-i}}\label{app:prop:FM-constructVis-withu}


We break the proof  into the following lemmas.

\begin{lemma}\label{lem:isotopy-2disk}
Fix an orientation of $S^2$. Let $\hat b_1, \hat b_2: D^2\to S^2$ be a pair of mutually disjoint, equi-oriented embeddings.  Let $N_1$, $N$ be 2-disk neighborhoods of $\hat b_1(D^2)$ and $\hat b_1(D^2)\cup \hat b_{2}(D^2)$ in $S^2$, respectively.
\begin{enumerate}[label=\textup{(\roman*)}, ref=\textup{(\roman*)},align=CenterWithParen]
\item There is an ambient isotopy $g:S^2\times I\to S^2$ with support on $N_1$ such that   $h_1\circ \hat b_1(x,y) = \hat b_1(-x,-y)$ for $(x,y)\in D^2$.
 \item There is an ambient isotopy $h:S^2\times I\to S^2$ with support on $N$\mycomment{i.e.,$h_t|_{\d N} = \id_{\d N}$ for all $t\in I$} such that $h_1\circ \hat b_1 = \hat b_2$ and $h_1\circ \hat b_2 = \hat b_1$.
\end{enumerate}
\end{lemma}

\begin{proof}
We prove (ii) only; (i) is easier (note that the transformation from $D^2=D^1\times D^1$ to itself given by $(x,y)\mapsto (-x,-y)$ is orientation-preserving).

Let $\hat N$ be a 2-disk neighborhood of $\hat b_1(D^2)\cup \hat b_{2}(D^2)$ in the interior of $N$, and choose a collar  $c:\d \hat N\times I\to N$ of $N$  such that  $c(x,0)=x$ for $x\in \d\hat N$ and $c(\d \hat N\times 1)=\d N$. Since the embeddings $\hat b_1$ and $\hat b_2$ are equi-oriented, by the Disk Theorem \cite[Corollary 3.3.7]{Kosinski} and the Isotopy Extension Theorem \cite[Theorem 2.5.2]{Kosinski}, there is an ambient isotopy $\hat h:\hat N\times I\to \hat N$ such that $\hat h_1\circ \hat b_1 = \hat \hat b_2$ and $\hat h_1\circ \hat b_2 = \hat b_1$. Choose a smooth function $m:I\to I$ satisfying $m(0)=1$ and $m(1)=0$, and define $h:S^2\times I\to S^2$  as follows. For each $t\in I$, let $h_t = \hat h_t$ on $\hat N$, let
\[
    h(c(x,s),t) = c(\hat h(x,m(s)t), s)
\]
for $(x,s)\in \d \hat N\times I$, and let $h_t=\id_{S^2}$ elsewhere.  It is readily verified that $h_0=\id_{S^2}$, and that for each $t\in I$,  the map $h_t$  is well-defined on $\d \hat N=c(d\hat N\times 0)$ and constant on the complement of $S^2\setminus \int N$.\mycomment{Note that $h$ is well-defined on $\d \hat N=c(\d \hat N\times 0)$, since for $s=0$ and $x\in \d\hat N$, $x = c(x,0)$  and $h(c(x,0),t)=c(\hat h(x,t),0) = c(x,0)=x = \hat h_t(x)$; $h$ is constant on the boundary, since when $s=1$ and $y\in \d N$ we have $y=c(x,1)$ for unique $x\in \d\hat N$, and so $h(y,t) = h(c(x,1),t) = c(\hat h(x,0), 1) = c(x,1) = y$. And when $t=0$, $h(c(x,s),0)=c(\hat h(x,m(s)t), s)=c(\hat h(x,0), s) = c(x,s)$.}%
\end{proof}

We extend these isotopies to  $S^4$ as follows.

\begin{lemma}\label{lem:isotopy-4disk}
Suppose that $b_1, b_2:D^4\to S^4$ are a pair of equi-oriented embeddings with mutually disjoint images such that, if $S^2\subset S^4$ denotes the standard embedding, we have $b_i^{-1}(S^2) = D^2\times0\times 0$ for $i=1,2$. Let $N_1$, $N$ be 2-disk neighborhoods of $b_1(D^4)\cap S^2$ and $(b_1(D^4)\cup b_2(D^4))\cap S^2$ in $S^2$, respectively.
\begin{enumerate}[label=\textup{(\roman*)}, ref=\textup{(\roman*)},align=CenterWithParen]
\item There is an ambient isotopy $F:S^4\times I\to S^4$ with support on an arbitrarily small 4-ball neighborhood of $b_1(D^4)\cup N_1$ such that $F_1$ fixes $S^2$ set-wise,  and $F_1(b_1(x,y)) = b_1(-x,y)$ for $(x,y)\in D^2\times D^2$.

\item There is an ambient isotopy $H:S^4\times I\to S^4$ with support on an arbitrarily small 4-ball neighborhood of $b_1(D^4)\cup b_2(D^4)\cup N$ such that $H_1$ fixes
\[
    S^2 \setminus \int \bigl(b_1(D^4)\cup b_2(D^4)\bigr) 
\]
set-wise, and  $H_1\circ b_1 = b_2$ and $H_1\circ b_2 = b_1$.
\end{enumerate}
\end{lemma}
\begin{proof}
We prove (ii) only; (i) is an analogous application of   part (i) of Lemma \ref{lem:isotopy-2disk}. Denote the closed 2-disk in $\R^2$ of radius $\tfrac{1}{2}$ by $\hat D^2$. Since for $i=1,2$, $b_i(D^4)$  intersects the standard 2-sphere along $b_i(D^2\times 0\times 0)$, we may identify a tubular neighborhood of $S^2=S^2\times 0\times 0$ with $S^2\times D^2$ so that there are equi-oriented, disjoint embeddings $\hat b_i:D^2\to b_i(D^4)\cap S^2$ such that  $b_i(D^4)=\hat b_i(D^2)\times \hat D^2$ and $b_i$ is given by $b_i(x,y)=(\hat b_i(x), \frac{1}{2}y)$ for $(x,y)\in D^2\times D^2$.

By Lemma \ref{lem:isotopy-2disk}(ii) there is an ambient isotopy $h:S^2\times I\to S^2$ with support on $N$ such that $h_1\circ \hat b_1 =  \hat b_2$ and $h_1\circ \hat b_2 = \hat b_1$; in particular, $h_1$ fixes $N\setminus \int \bigl(\hat b_1(D^2)\cup \hat b_2(D^2)\bigr)$ set-wise. We construct an isotopy $H$ of $N\times D^2\subset S^2\times D^2$ such that, for each $t\in I$:
\begin{itemize}
    \item[(1)] $H_t|_{\d(N\times D^2)} = \id_{\d (N\times D^2)}$,
    \item[(2)] $H_t(x,y) = (h_t(x),y)$ \mycomment{$H(x,y,t) = (h(x,t),y)$}  for all $x\in N$, $y\in \hat D^2$, and\mycomment{ H|N\times \hat D^2\times I = h|N\times I \times id|\hat D^2}
    \item[(3)] $H_1\circ b_1 = b_2$ and $H_1\circ b_2 = b_1$.
\end{itemize}
Choose a smooth function $m:I\to I$ such that $m(1)=0$ and $m(s)=1$ for all $s\in[\tfrac{1}{2},1]$. For each $t\in I$  and $(x,y)\in N\times D^2$, let $H_t(x,y) = (h_{m(|y|)t}(x), y)$, where $|\cdot |$ denotes the Euclidean norm on $D^2$. Note that on $N\times D^2$, $H_t$  has inverse given by ${H_t}^{-1}(x,y) = (h_{m(|y|)t}^{-1}(x), y)$. To verify (1), observe that for $x\in \d N$ and $y\in D^2$ we have $(h_{m(|y|)t}(x), y)=(x, y)$; for $x\in N$ and $y\in \d D^2$, we have $(h_{m(|y|)t}(x), y)=(h_{0}(x), y)=(x,y)$. To verify (2), observe that if $y\in \hat D^2$ then $m(|y|)=1$ and so $(h_{m(|y|)t}(x), y) = (h_t(x), y)$. Regarding (3), for  $(x,y)\in D^2\times D^2$ we have
\[
    H_1(b_1(x,y)) = H_1(\hat b_1(x), \tfrac{1}{2}y) = (h_1(\hat b_1(x)),\tfrac{1}{2}y) = (\hat b_2(x),\tfrac{1}{2}y) = b_2(x,y),
\]
and we have $H_1\circ b_2=b_1$ similarly.

Now, by (1) we may extend $H$ to an isotopy of $S^4$ that is constant on the complement of $N\times D^2$.  Since $h_1$ fixes $N\setminus \int \bigl(\hat b_1(D^2)\cup \hat b_2(D^2)\bigr)=\bigl[N\setminus \int \bigl(\hat b_1(D^2)\cup \hat b_2(D^2)\bigr)\bigr]\times 0$ set-wise, so does $H_1$  by property (2); hence $H_1$ fixes $S^2$ (since $H_1$ is the identity outside $N^2\times D^2$).%
\end{proof}

We may now prove Lemma \ref{prop:exchange-all-b-i}.

\begin{proof}[Proof of Lemma \ref{prop:exchange-all-b-i}] \mycomment{ disk theorem => unknotted unique}
Choose an orientation-preserving diffeomorphism $\Phi:S^4\to S^4$ that takes $\U$ to the standard embedding $S^2\subset S^4$; then $b_i'=\Phi\circ b_i$, $i=1,\ldots, d$, is a collection of mutually disjoint, equi-oriented embeddings $D^4\to S^4$ whose images intersect $S^2$ precisely along $b_i'(D^2\times 0\times 0)$, respectively.

If $\rho$ is non-trivial, write it as a product of non-trivial transpositions $\rho=\tau_1\tau_2\ldots \tau_n$ for some $n\geq 1$. For each $1\leq k\leq n$, write $\tau_k=(a_k\; b_k)$ for some $a_k,b_k\in \{1,2,\ldots, d\}$, and let $N\sup{k}$ be a 2-disk neighborhood of $(b_{a_k}'(D^4)\cup b_{b_k}'(D^4))\cap S^2$ in $S^2\setminus \int \mycup{i\neq a_k, b_k}{} b_i'(D^4)$. By Lemma \ref{lem:isotopy-4disk}(ii) there is an ambient isotopy $H^{(k)}:S^4\times I\to S^4$ with support on $N\sup{k}\times D^2$ in $S^4\setminus \int \mycup{i\neq a_k, b_k}{} b_i'(D^4)$ such that $H\sup{k}_1$ fixes
\[
    S^2 \setminus \int \bigl(b_{a_k}(D^4)\cup b_{b_k}(D^4)\bigr) 
\]
set-wise and is such that $H\sup{k}_1\circ b_{i}' = b_{\tau_k(i)}'$ for $i=a_k, b_k$.  Define $H:S^4\times I\to S^4$ by
\[
    H(x,t)=H\sup{k}(x,{n(t-\tfrac{k-1}{n})})
\]
for $x\in S^4$ and $t\in [\tfrac{k-1}{n}, \tfrac{k}{n}]$, where $k=1,2,\ldots, n$. Then $H$ is an ambient isotopy which fixes
\[
    S^2\setminus \int \mycup{i=1}{d} b_i'(D^4)
\]
set-wise and satisfies
\[
   H_1\circ b_i' = b_{\rho(i)}'
\]
for each $1\leq i\leq d$. Now, by Lemma \ref{lem:isotopy-4disk}(i), for each $1\leq i\leq d$ there is an ambient isotopy $F^{(i)}:S^4\times I\to S^4$ with support on a 4-ball neighborhood of $N\sup{\rho(i)}$ in $S^4\setminus \int \mycup{k\neq \rho(i)}{} b_k'(D^4)$ such that $F\sup{i}_1$ fixes $S^2$ set-wise and is such that
\[
    F\sup{i}_1\circ b_{\rho(i)}'(x,y) = b_{\rho(i)}'(\mu_i \, x,y)
\]
for $(x,y)\in D^2\times D^2$.  Define $F:S^4\times I\to S^4$ by
\[
    F(x,t)=F\sup{i}(x,{d(t-\tfrac{i-1}{d})})
\]
for $x\in S^4$ and $t\in [\tfrac{i-1}{d}, \tfrac{i}{d}]$, where $i=1,2,\ldots, d$. Then $F$ is an ambient isotopy which fixes $S^2$ set-wise and satisfies
\[
    F_1\circ b_{\rho(i)}'(x,y) = b_{\rho(i)}'(\mu_i \, x,y)
\]
for $(x,y)\in D^2\times D^2$. Thus if $K:S^4\times I\to S^4$ is the ambient isotopy defined by $H_{2t}$ for $t\in [0,\tfrac{1}{2}]$ and $F_{2t-1}$ for $t\in [\tfrac{1}{2},1]$, then $\hat\varphi_t = \Phi^{-1}\circ K_t\circ \Phi$ is the required isotopy.
\end{proof}
%


\begin{thebibliography}{test}

%
\bibitem{B} A. Bartels,  Higher dimensional links are singular slice, \textit{Math. Annal.} \textbf{320} (3) (2001) 547-576.
%
\bibitem{BT} A. Bartels and P. Teichner, All two dimensional links are null homotopic, \textit{Geom. Topol.} \textbf{3} (1999) 235-252.
%
%
%
%
\bibitem{FR} R. Fenn and D. Rolfsen,  Spheres may link in 4-space, \textit{J. London Math. Soc.} \textbf{34} (1986) 177-184.
%
\bibitem{FK} M. Freedman and R. Kirby,  A  geometric  proof  of  Rochlin's  theorem, \textit{Proc. Symp.  Pure
Math. AMS}  \textbf{32}  (2) (1978)  85-98.
%
\bibitem{FQ} M. Freedman and F. Quinn,  The topology of 4-manifolds, \textit{Princeton Math. Series,} Vol. 39 (Princeton, NJ, 1990).
%
%
\bibitem{HL} N. Habegger and X.-S. Lin,  The classification  of  links up to homotopy, \textit{J.  Amer.  Math. Soc.} \textbf{3} (1990) 389-420.
%
%
%
\bibitem{Hos} F. Hosokawa and A. Kawauchi,  Proposals for unknotted surfaces in four-spaces, \textit{Osaka J. Math.} \textbf{16} (1979) 233-248.
%
\bibitem{K} S. Kamada, Vanishing of a certain kind of Vassiliev invariants of 2-knots, \textit{Proc. AMS}
 \textbf{127} (11) (1999) 3421-3426.
%
%
\bibitem{Ki1}  P. Kirk,  Link maps in the four sphere, in \textit{Proc. 1987 Siegen Topology Conf.}, SLNM 1350 (Springer, Berlin, 1988).
%
%
\bibitem{Ki2} P. Kirk,  Link homotopy with one codimension two component. \textit{Trans. Amer. Math.
Sot.} \textbf{319} (1990) 663-688.
%

\bibitem{Ko1}  U. Koschorke, On link maps  and  their  homotopy  classification, \textit{Math.  Ann.} \textbf{286} (1990) 753-782.
%
\bibitem{Kosinski} A. Kosinski,  \textit{Differential manifolds} (Academic Press, San Diego, 1993).
%
%
\bibitem{Li97}  G.-S. Li,  An invariant of link homotopy in dimension four, \textit{Topology} \textbf{36} (1997) 881-897.
%
%
%
\bibitem{me2} A. Lightfoot,  On invariants of link maps in dimension four, \textit{ J. Knot Theory Ramifications,} July 2016, \\DOI: http://dx.doi.org/10.1142/S0218216516500607.
%
%
\bibitem{Mi} J. Milnor,  Link groups, \textit{Ann. of Math} \textbf{59} (1954) 177-195.
%
\bibitem{MR}  W. S. Massey and D. Rolfsen,  Homotopy classification of higher-dimensional links, \textit{Indiana Univ. Math. J.}  \textbf{34} (1985)  375-391.

%
\bibitem{Pilz}  A. Pilz,  Verschlingungshomotopie von 2-Sph\"{a}ren im 4-dimensionalen Raum, Diploma thesis,  University of Seigen (1997).

\bibitem{Sc}  G. P. Scott,  Homotopy links, \textit{Abh. Math. Sem. Univ. Hamburg} \textbf{32} (1968) 186-190.

%
\bibitem{ST} R. Schneiderman and P. Teichner,  Higher order intersection numbers of 2-spheres in
4-manifolds, \textit{Algebraic and Geometric Topology} \textbf{1} (2001) 1-29.
%
%
%
\bibitem{W}  C. T. C. Wall,  \textit{Surgery on Compact Manifolds} (Academic Press, New York, 1970).

\end{thebibliography}
\addcontentsline{toc}{chapter}{Bibliography}

\end{document}